\documentclass[11pt]{article}
\usepackage{amsmath}
\usepackage{amssymb}
\usepackage{ntheorem}
\usepackage{graphicx}
\usepackage{subfigure}
\usepackage{color}
\usepackage{bm}
\usepackage{comment}

\numberwithin{equation}{section}
\newtheorem{theorem}{Theorem}[section]
\newtheorem{definition}[theorem]{Definition}
\newtheorem{lemma}[theorem]{Lemma}
\newtheorem{corollary}[theorem]{Corrolary}
\newtheorem{remark}[theorem]{Remark}
\def\Proof{\noindent{\bf Proof} \quad}
\def\qed{\hfill $\Box$ \smallskip}

\def\Z{{\mathbb Z}}

\def\R{{\mathbb R}}

\renewcommand{\l}{\left}
\renewcommand{\r}{\right}

\renewcommand{\d}{\mathrm{d}}
\newcommand{\be}{\begin{equation}}
	\newcommand{\ee}{\end{equation}}

\newcommand\cN{{\mathcal N}}
\renewcommand{\Z}{{\mathbb Z}}
\newcommand\cG{{\mathcal G}}

\renewcommand{\S}{{\mathbb S}}
\renewcommand{\R}{{\mathbb R}}
\newcommand\cT{{\mathcal T}}

\newcommand{\p}{\partial}
\renewcommand{\d}{\mathrm{d}}
\newcommand{\eps}{\varepsilon}

\begin{document}

\title{Convergence analysis of three semi-discrete numerical schemes for nonlocal geometric flows including perimeter terms}

\author{Wei, Jiang\thanks{School of Mathematics and Statistics, Wuhan University, Wuhan, 430072, China ({jiangwei1007@whu.edu.cn}).} \and Chunmei, Su \thanks{Yau Mathematical Sciences Center, Tsinghua Unversity, Beijing, 100084, China({sucm@tsinghua.edu.cn}).}\and Ganghui, Zhang \thanks{Yau Mathematical Sciences Center, Tsinghua Unversity, Beijing, 100084, China({gh-zhang19@mails.tsinghua.edu.cn}).}}

\maketitle

\begin{abstract}
   We present and analyze three distinct semi-discrete schemes for solving  nonlocal geometric flows incorporating perimeter terms. These schemes are based on the finite difference method, the finite element method, and the finite element method with a specific tangential motion. We offer rigorous proofs of quadratic  convergence under $H^1$-norm for the first scheme and linear convergence under $H^1$-norm for the latter two  schemes. All error estimates rely on the observation that the error of the nonlocal term can be controlled by the error of the local term. Furthermore, we explore the relationship between the convergence under $L^\infty$-norm and manifold distance.  Extensive numerical experiments are conducted to verify the convergence analysis, and demonstrate  the accuracy of our schemes under various norms for different types of nonlocal flows.
\end{abstract}

{\small {\bf MSCcodes}: 
65M60, 65M12, 35K55

}


\section{Introduction} 

In this paper, we analyze and establish the convergence result of three distinct numerical methods for evolving a closed plane curve $\Gamma(t)$ under a nonlocal flow that involves perimeter. The normal velocity of $\Gamma(t)$ is determined by the formula
\begin{equation}\label{Nonlocal}
\mathcal{V}= \l(\kappa-f(L)\r)\cN, \end{equation}
where $\kappa$ represents the curvature of the curve, $f$ is a Lipschitz function, $L$ is the perimeter, and $\cN$ is the unit inner normal vector. Equation \eqref{Nonlocal} encompasses a wide range of geometric flows, including:
\[f(L)= \begin{cases}
\frac{2\pi}{L}, &\text{for area-preserving curve shortening flow of simple curves \cite{Gage1986} },\\
 \frac{2\pi\, \mathrm{ind}(\Gamma)}{L}, &\text{for area-preserving curve shortening flow of nonsimple curves \cite{Wang2014}}, \\
  \frac{2\pi-\beta}{L}, &\text{for curve flows with a prescribed rate of change of area \cite{Dallaston2016,Tsai2018}},  \end{cases} \]
where $\beta\in(-\infty,\infty)$, $\mathrm{ind}(\Gamma)\in \Z$ denotes the rational index \cite{Carmo2016} of a nonsimple curve $\Gamma$. The inclusion of an additional nonlocal force,
 $f(L)$, enables us to control the area change of an evolving curve. Indeed, by the theorem of turning tangents \cite{Carmo2016}, the rate of area change can be determined by \cite{Deckelnick2005}
\[
\frac{\d A}{\d t}=\int_{\Gamma}\mathcal{V}\cdot \cN \d s=\begin{cases}
    0,\quad &\text{for}\quad f(L)=\frac{2\pi }{L}\quad \text{and simple curves} ,\\
	0,\quad &\text{for} \quad f(L)=\frac{2\pi \,\mathrm{ind}(\Gamma)}{L}\quad \text{and nonsimple curves},\\
	-\beta,\quad &\text{for}\quad f(L)=\frac{2\pi-\beta}{L}\quad \text{and simple curves}.
\end{cases}
\]
In this paper, we focus on the study of curve evolutions that maintain their topological characteristics.

In recent years, there has been significant emphasis on the development of theoretical and modeling frameworks for nonlocal geometric flows. One prominent instance of such work is the area-preserving curve shortening flow (AP-CSF), which has become vital in the field of image processing \cite{Dolcetta2002,Sapiro1995,Sapiro2001} and can be interpreted as a limit of the nonlocal model of the Ginzburg-Laudau equation \cite{Bronsard1997}.  The existence and convergence results of AP-CSF for both simple and nonsimple closed curve cases have been extensively explored \cite{Gage1986,Wang2014}. Moreover, the study of curve flows with a prescribed rate of change in the enclosed area has arisen in connection with the investigation of contracting bubbles in fluid mechanics \cite{Dallaston2012,Dallaston2013,Dallaston2016}, and the long-time behavior of such flows has been addressed in \cite{Tsai2018}. For more comprehensive theoretical studies related with nonlocal geometric flows, we refer to \cite{Chambolle2015,Jiang2008,Tsai2015}.

Extensive numerical methods have been employed to simulate the AP-CSF and curve flows with a prescribed rate of change in the enclosed area. Examples of such methods for the AP-CSF include the finite difference method \cite{Mayer2000}, the MBO method \cite{Ruuth2003}, the crystalline algorithm \cite{Ushijima2004}, as well as PFEMs  \cite{Pei2022,Barrett2020}. Additionally, a rescaled spectral collocation scheme was proposed in \cite{Dallaston2016} for closed embedded plane curves
with a prescribed rate of change in the enclosed area. However, there has been relatively little research on the numerical analysis of these methods. Recently, in \cite{Jiang2022}, the authors proposed a semi-discrete finite element method for the AP-CSF of simple curves and established its convergence in $H^1$-norm. The nonlocal nature of the geometric equations presents a major challenge for the error analysis.

In this paper, we propose three numerical schemes for  nonlocal geometric flows involving perimeter \eqref{Nonlocal} and give their error estimates. Our main observation is that the difference between the nonlocal term  and its discrete version  can be managed through the disparity of the local term. Specifically, we introduce the following three different types of semi-discrete  schemes:
\begin{itemize}
	\item Firstly, we employ a finite difference method to discretize the parametrization equation of \eqref{Nonlocal}
	\begin{equation}\label{Nonlocal,parametrization1}
\begin{split}
&\p_t X=\frac{1}{|\p_\xi X|}\p_\xi\Big(\frac{\p_\xi X}{|\p_\xi X|} \Big)-f(L)\Big(\frac{\p_\xi X}{|\p_\xi X|}\Big)^{\perp},\quad \xi\in \S^1,
\end{split}
 \end{equation}
where $\S^1=[0, 2\pi]$, $(a,b)^\perp:=(-b,a)$ denotes a clockwise rotation by $\pi/2$ and the periodic function $X(\xi, t): \S^1\times[0,T] \rightarrow \mathbb{R}^2$ is a parameterization of the closed curve $\Gamma \subset \mathbb{R}^2$. Under certain appropriate assumptions, we demonstrate that  the resulting semi-discrete scheme converges quadratically in the discrete $H^1$-norm as defined in \cite{Deckelnick2022}. The proof is based on a careful Taylor expansion result and an averaged approximation of the normal vector.
\smallskip

\item Secondly, we utilize a finite element method for a natural weak formulation of  \eqref{Nonlocal,parametrization1}. The derived semi-discrete scheme is based on our previous work on AP-CSF of simple curves \cite{Jiang2022}. An $H^1$-optimal  error estimate follows from our key observation mentioned above.
    \smallskip

\item Thirdly, we introduce an artificial  tangential motion and apply a finite element method for an alternative parametrization of the geometric equation
\begin{equation}\label{Nonlocal1}
	\p_t X = \frac{\p_{\xi\xi}X}{|\p_\xi X |^2}  -f(L)\cN.
\end{equation}
This form of reparametrization was initially  proposed by Deckelnick and Dziuk for the curve shortening flow \cite{Deckelnick1995} to improve the mesh quality during evolution. It was later interpreted as a DeTurck trick by Elliot and Fritz in \cite{Elliott2017}. Recently, the DeTurck trick has been further applied to various geometric flows such as elastic flow  \cite{Pozzi2022}, anisotropic curve shortening flow \cite{Deckelnick2023,Deckelnick2023b,Deckelnick2023c} and fourth-order flows \cite{Deckelnick2024}. We emphasize that we have successfully extended the DeTurck trick to the general nonlocal flow case. The resulting semi-discrete scheme yields an asymptotic equidistribution property, as well as an $H^1$-optimal error estimate.
\end{itemize}

As a byporduct, we further explore the convergence of the schemes under manifold distance,  a topic extensively discussed in the numerical computation community \cite{Bao2023,Bao2021,Zhao2021,Jiang2023,Jiang2024stable}. We prove that, for simple curves, convergence in the function $L^\infty$-norm implies convergence under the manifold distance. Moreover, we prove an optimal convergence of the finite difference scheme under the manifold distance.

The rest of this paper is organized as follows.  In Section 2, we briefly introduce the mathematical notations. In Section 3, we propose the semi-discrete schemes and provide the error estimates  for the  finite difference method. In Section 4, we consider  the finite element method, and the finite element method with a tangential motion.  Section 5 aims to establishing a connection between the convergence of the manifold distance and $L^\infty$-norm. Section 6 presents extensive numerical experiments for the three different numerical schemes and various types of nonlocal flows. The numerical results demonstrate our convergence analysis results in both the $H^1$-norm and the manifold distance. Moreover, a better mesh quality is achieved for the finite element method with the aid of tangential motions. Finally, we draw some conclusions in Section 7.

\section{Notations}

Throughout the paper, we denote the quantities related to the real solution and discrete solution by capital and lower-case letters, respectively. Specifically, for the solution of \eqref{Nonlocal,parametrization1},
we denote $\cT=\frac{\p_\xi X}{|\p_\xi X|}$ and $\cN=\cT^\perp$ by the tangent and inner normal of the curve, respectively. Thus \eqref{Nonlocal,parametrization1} can be simply written as
\begin{equation}\label{para2}
\p_t X=\frac{1}{|\p_\xi X|}\p_\xi\cT -f(L)\cN,\quad \xi\in \S^1;\qquad
X(\xi,0)=X^0(\xi).
		 \end{equation}
Direct computation gives
\be\label{leg}
\p_t|\p_\xi X |=\p_{t}\p_\xi X\cdot \cT=\p_\xi(\p_{t} X\cdot \cT)-\p_{t} X\cdot \p_\xi\cT= -|\p_\xi X||\p_t X |^2-f(L)|\p_\xi X|\p_t X\cdot \mathcal{N}.
\ee

For spatial discretization, we utilize a uniform mesh, where the equidistributed grid points $\cG_h:=\{\xi_1,\ldots,\xi_N\}\subset\S^1$ are given by $\xi_j=jh,j=0,\ldots,N$ for $h=2\pi/N$ with $N\ge 2$. We use a periodic index, i.e., $a_{j}=a_{j\pm N}$ when involved. Denote $X_j=X(\xi_j)$,  $\dot{X}_j=\partial_t X(\xi_j)$, and set
\[Q_j=|X_{j}-X_{j-1}|, \quad \mathcal{T}_{j}=\frac{X_{j}-X_{j-1}}{Q_j} ,\quad j=1,\ldots,N.
\]
Let $x_h:\cG_h\rightarrow \R^2$ be a grid function. We define the discrete length element $q_j$, the discrete tangent $\tau_j$ and normal $n_j$ as
\be\label{qtau}
q_j=|x_j-x_{j-1} | ,\quad \tau_j=\frac{ x_j-x_{j-1}}{q_j},\quad n_j=\tau^\perp_j,
\ee
where $x_j=x_h(\xi_j)$ denotes the vertex of the polygon that approximates the curve. Denote $l_h=\sum\limits_{j=1}^N q_j$ by the perimeter of the polygon.
Throughout the article, we denote $C$ by a general constant which is independent of the mesh size $h$ and might vary from line to line.

\smallskip

\textbf{(Assumption 2.1)} Suppose that the solution of \eqref{Nonlocal,parametrization1} satisfies $X\in C^1\l([0,T],C^4(\S^1)\r)$, i.e., \[K_1(X):=\|X\|_{C^1\l([0,T],C^4(\S^1)\r)}<\infty,\]
	and there exist constants  $0<C_1<C_2$  such that
\be\label{C12}
C_1\le \l|\p_\xi X(\xi,t)\r|\le C_2,\quad \forall\ (\xi, t)\in \S^1\times [0, T].
	  \ee
Under this assumption, we have the following results, which have been established in \cite{Deckelnick2022}.

\begin{lemma}\cite[Lemmas 3.1, 3.3]{Deckelnick2022}
Under Assumption 2.1, there exists $h_0>0$ such that for  $0<h\le h_0$, we have
\begin{subequations}\label{Taylor expansion}
			\be\label{Qes}
				C_1 \leq Q_j/h \leq C_2,\quad
	\frac{Q_j+Q_{j+1}}{2h}  = | \p_\xi X(\xi_j) | + O(h^2),
\ee
\be\label{cTe}
\cT_{j} + \cT_{j+1}= 2 \, \cT(\xi_j) + O(h^2),\quad \cT_j =\frac{1}{2}\l(\cT(\xi_j)+\cT(\xi_{j-1}) \r)+O(h^2), \ee
\be\label{Xes}
	\frac{\cT_{j+1} - \cT_j}{h}  = \p_\xi\cT(\xi_j) + O(h^2),	\quad \frac{\dot X_{j+1}-\dot X_{j}}{h}=\frac{1}{2}\l(\p_{t}\p_\xi X(\xi_{j})+\p_{t}\p_\xi X(\xi_{j+1})\r)+ O(h^2),
			\ee
\be\label{tausT}
			\l|\tau_{j+1/2}^\perp -\cN(\xi_j)\r|\le \frac{2}{|\cT_{j}+\cT_{j+1}|}\l(|\cT_{j}-\tau_{j}|+|\cT_{j+1}-\tau_{j+1}|\r)+Ch^2,
			\ee
		\end{subequations}
where $\tau_{j+1/2}:=\frac{\tau_j+\tau_{j+1}}{|\tau_j+\tau_{j+1}|}$ represents the averaged vertex tangent.
\end{lemma}

For a grid function $u:\cG_h\rightarrow \R^2$, we define the backward  difference quotient as
\[
\delta u_j :=\frac{u_j-u_{j-1}}{h},\quad j=1,\ldots,N.
\]
Moreover, to measure the error, we introduce the following discrete norms:
\begin{equation}\label{Definition of grid norm}
	\|u\|_{L^2_G}: = \Bigl( h \sum_{j=1}^N |u_j |^2 \Bigr)^{\frac{1}{2}}, \quad
\| u\|_{H^1_G}:= \Bigl( h \sum_{j=1}^N \bigl( | u_j |^2 + | \delta u_j |^2 \bigr) \Bigr)^{\frac{1}{2}}.
\end{equation}

\section{Finite difference method}\label{sec2}
 In this section, we utilize a finite difference method to solve the equation \eqref{Nonlocal,parametrization1}.
\begin{definition}
A semi-discrete finite difference approximation of \eqref{Nonlocal,parametrization1} is to find a grid function $x_h: \cG_h\times [0,T]\rightarrow \R^2$ such that	
\begin{equation}\label{Semidiscrete,FDM}
\dot{x}_j = \frac{2}{q_j + q_{j+1}}
\l(\tau_{j+1} - \tau_{j}\r) -f(l_h)\tau^\perp_{j+1/2}
\quad \text{in }\quad  (0,T];\qquad
x_j(0) = X^0(\xi_j).
\end{equation}
\end{definition}

\begin{theorem}\label{Error estimate 3}
Let $X(\xi,t)$ be a solution of \eqref{Nonlocal,parametrization1}  satisfying Assumption 2.1. Then there exists $h_0>0$ such that for all $0<h\le h_0$, there exists a unique finite difference semi-discrete solution $x_h(t)$  in the sense of \eqref{Semidiscrete,FDM}. Furthermore, we have the following  error estimate
\begin{equation}\label{Main estimate 3}
\sup_{t\in [0,T]} \| X(t) - x_h(t)  \|_{H^1_G}   \leq C h^2,
\end{equation}
where  the constants  $h_0,C$ depend on $C_1,C_2,K_1(X),T$ and $f$.
\end{theorem}

 First, we compute the evolution equation for the discrete length $q_j$.

\begin{lemma}
	Suppose $x_h$ is the finite difference semi-discrete solution  in the sense of   \eqref{Semidiscrete,FDM}, then we have
\begin{equation}\label{FDM, length equation}
	\dot{q}_j
+ \tfrac{q_{j-1} + q_j}{4} | \dot{x}_{j-1} |^2
+ \tfrac{q_j + q_{j+1}}{4}| \dot{x}_j |^2+\tfrac{q_{j-1}+ q_j}{4}f(l_h)\dot x_{j-1}\cdot \tau_{j-1/2}^\perp +\tfrac{q_{j+1}+ q_j}{4}f(l_h)\dot x_{j}\cdot \tau_{j+1/2}^\perp= 0.
\end{equation}
\end{lemma}

\Proof
~We begin by computing $\dot{q}_j $ as
\begin{align}
	\dot{q}_j
	&=\l(\dot x_j-\dot x_{j-1}\r)\cdot \tau_j\notag\\
	&=\tau_j\cdot \l(\frac{2}{q_j + q_{j+1}}
\l(\tau_{j+1} - \tau_{j}\r) -f(l_h)\tau_{j+1/2}^\perp-\frac{2}{q_j + q_{j-1}}
\l(\tau_{j} - \tau_{j-1}\r) +f(l_h)\tau_{j-1/2}^\perp \r)\notag\\
&=\frac{2}{q_j + q_{j+1}}\l(\tau_j\cdot \tau_{j+1} -1\r)-f(l_h)\tau_{j+1/2}^\perp\cdot \tau_j
 -\frac{2}{q_j + q_{j-1}}\l(1-\tau_j\cdot \tau_{j-1}\r)+f(l_h)\tau_{j-1/2}^\perp\cdot \tau_j\notag\\
 &=\frac{2}{q_j + q_{j+1}}\l(\tau_j\cdot \tau_{j+1} -1\r)-f(l_h)\tau_{j+1/2}^\perp\cdot \tau_j
 +\frac{2}{q_j + q_{j-1}}\l(\tau_j\cdot \tau_{j-1}-1\r)-f(l_h)\tau_{j-1/2}^\perp\cdot \tau_{j-1}\notag\\
 &=:J_j+J_{j-1},\label{qjpf}
\end{align}
where for the last second equality, we have employed the property
\begin{equation}\label{taujh}
	\tau_{j-1/2}^\perp\cdot \tau_{j}
	=\tau_{j-1}^\perp/|\tau_j+\tau_{j-1}| \cdot \tau_j=-\tau_{j-1}\cdot \tau_j^\perp/|\tau_j+\tau_{j-1}|=-\tau_{j-1}\cdot\tau_{j-1/2}^\perp.
\end{equation}
Multiplying  \eqref{Semidiscrete,FDM} by $ \frac{q_j + q_{j+1}}{4}f(l_h)\tau_{j+1/2}^\perp$, we obtain
\[
	\frac{q_j + q_{j+1}}{4} f(l_h)\tau_{j+1/2}^\perp\cdot \dot x_j+\frac{q_j + q_{j+1}}{4}f(l_h)^2-f(l_h)\tau_{j+1/2}^\perp\cdot \frac{\tau_{j+1}-\tau_j}{2}=0,
\]
which can be simplified as
\begin{equation}\label{FDM,length2}
	\frac{q_j + q_{j+1}}{4} f(l_h)\tau_{j+1/2}^\perp\cdot \dot x_j+\frac{q_j + q_{j+1}}{4}f(l_h)^2+f(l_h)\tau_{j+1/2}^\perp\cdot \tau_j=0,
\end{equation}
by using \eqref{taujh}. Combining  \eqref{Semidiscrete,FDM} and \eqref{FDM,length2}, we get
\begin{align*}
	J_j
	&=\frac{2}{q_j + q_{j+1}}\l(-\frac{1}{2}|\tau_j-\tau_{j+1}|^2\r)-f(l_h)\tau_{j+1/2}^\perp\cdot \tau_j\\
	&=-\frac{1}{q_j + q_{j+1}}\l|\dot x_j+f(l_h)\tau_{j+1/2}^\perp \r|^2\l(\frac{q_j + q_{j+1}}{2} \r)^2-f(l_h)\tau_{j+1/2}^\perp\cdot \tau_j\\
	&=-\frac{q_j + q_{j+1}}{4}\l(|\dot x_j|^2+2 f(l_h)\tau_{j+1/2}^\perp\cdot \dot x_j+f(l_h)^2 \r)-f(L_h)\tau_{j+1/2}^\perp\cdot \tau_j\\
	&=-\frac{q_j + q_{j+1}}{4}|\dot x_j|^2-\frac{q_j + q_{j+1}}{2} f(l_h)\tau_{j+1/2}^\perp\cdot \dot{x}_j-\frac{q_j + q_{j+1}}{4}f(l_h)^2-f(l_h)\tau_{j+1/2}^\perp\cdot \tau_j\\
	&=-\frac{q_j + q_{j+1}}{4}|\dot x_j|^2-\frac{q_j + q_{j+1}}{4} f(l_h)\tau_{j+1/2}^\perp\cdot \dot x_j.
\end{align*}
 Plugging this into \eqref{qjpf} yields \eqref{FDM, length equation}, and the proof is completed.
 \qed

\medskip

 \Proof[Proof of Theorem \ref{Error estimate 3}]
~We define
\begin{equation}\label{FDM,induction}
	\begin{split}
		T^* = \sup \Bigl\{ t \in [0,T]: x_h \,\,\,\mathrm{solves}\,\,\, \eqref{Semidiscrete,FDM} \,\,\,\text{ with }\,\,\, \tfrac{C_1}{2} \leq \frac{q_j(t)}{h} \leq 2 C_2, \,\,\,\max_{j=1,\ldots,N }
| \cT_j(t)-\tau_j(t) | \leq h^{\frac54}  \Bigr\}.
	\end{split}
\end{equation}
Clearly $ T^*>0$. Noticing the nonlinear terms in \eqref{Semidiscrete,FDM} are locally Lipschitz with respect to $x_j$, we get local existence and uniqueness using standard ODE theory. Furthermore, since $q_j(0)=Q_j(0)$ and $\tau_j(0)=\cT_j(0)$, the desired estimate also holds by continuity. By \eqref{FDM,induction} and the Lipschitz continuity of $f$, we have $\forall\ t\in [0,T^*]$,
 \begin{equation}\label{Control of perimeter,FDM}
 \begin{split}
 	2\pi C_1\le  L\le 2\pi C_2 ,&\quad \pi C_1\le l_h\le 4\pi C_2 ,\quad  |f(L)|\le C,
 \end{split}
 \end{equation}
 where $C$ is a constant depending on $C_1, C_2$ and $f$. We claim that there exists a constant $h_1>0$ such that for $0< h\le h_1$, it holds
\begin{equation}\label{Induction estimate 1}
	\max_{j=1,\ldots,N} | \dot{x}_j(t) |
\leq C,\quad \forall\ t\in [0,T^*],
\end{equation}
where $C$ depends on $C_1, C_2$ and $f$. Indeed, by \eqref{Semidiscrete,FDM}, \eqref{Xes}, \eqref{FDM,induction} and \eqref{Control of perimeter,FDM}, we obtain
\begin{align*}
	 |\dot{x}_j|^2
	 & \le  2\l|\frac{2}{q_j + q_{j+1}}
\l(\tau_{j+1} - \tau_{j}\r)\r|^2 +2\l|f(l_h)\tau_{j+1/2}^\perp \r|^2\le C\l|\frac{\tau_{j+1} - \tau_{j}}{h}\r|^2+C\\
&\le C\l(\l|\frac{\cT_{j+1} - \cT_{j}}{h}\r|+\frac{2}{h}\max_{k=1,\ldots,N}|\cT_{k}-\tau_{k} |\r)^2+C\le C.
\end{align*}
Moreover, based on \eqref{cTe}, we have
\begin{equation}\label{Induction estimate 2}
	\min_{j=1,\ldots,N}|\cT_{j}+\cT_{j+1}|\ge 1,
\end{equation}
when $h$ is sufficiently small.
Define the truncation error as
\begin{align}
	\mathcal{R}_j &:=  \dot{X}_j- \frac{2}{Q_j + Q_{j+1}}  \l( \mathcal{T}_{j+1} - \mathcal{T}_{j}\r)+f(L)\cN (\xi_j), \label{Remainder term1,AP-CSF} \\
\widetilde{\mathcal{R}}_j & := \dot{Q}_j  + \tfrac{ Q_{j-1}+ Q_j}{4}| \dot{X}_{j-1} |^2 +\tfrac{ Q_{j}+ Q_{j+1}}{4}| \dot{X}_j |^2 \notag\\
&\quad +\tfrac{Q_{j-1}+ Q_j}{4} f(L)\dot X_{j-1}\cdot \cN(\xi_{j-1}) +\tfrac{Q_{j+1}+ Q_j}{4}f(L)\dot X_{j}\cdot \cN(\xi_{j}).\label{Remainder term2,AP-CSF}
\end{align}

\begin{enumerate}
	\item[(1).] \emph{Estimates of the truncation error $\mathcal{R}_j,\widetilde{\mathcal{R}}_j$}.
Employing \eqref{Nonlocal,parametrization1}, \eqref{Qes} and \eqref{Xes}, one gets
	\begin{equation}\label{Estimate of Remainder 1}
		\begin{split}
			\mathcal{R}_j&= \dot{X}_j- \frac{2}{Q_j + Q_{j+1}}\l(\cT_{j+1} - \cT_{j}\r)+f(L)\cN(\xi_j)\\
			&=\dot{X}_j-\frac{1}{|\p_\xi X(\xi_j)| + O(h^2)}\l( \p_\xi\cT(\xi_j) + O(h^2)\r)+f(L)\cN(\xi_j)\\
			&=\dot{X}_j-\frac{1}{|\p_\xi X(\xi_j) |}\l(1+O(h^2)\r)\cdot \l( \p_\xi\cT(\xi_j) + O(h^2)\r)+f(L)\cN(\xi_j)\\
			&=\dot{X}_j-\frac{1}{|\p_\xi X(\xi_j) |}\p_\xi\cT(\xi_j)+f(L)\cN(\xi_j)+O(h^2)\\
			&= O(h^2).
		\end{split}
	\end{equation}
	Similarly, applying \eqref{leg}, \eqref{cTe} and \eqref{Xes}, we derive
		\begin{align*}
			\dot{Q}_j
			&=\l(\dot X_j-\dot X_{j-1}\r)\cdot \cT_j\\
			&=\frac{h}{4}\l(\p_{t}\p_\xi X(\xi_{j-1})+\p_{t}\p_\xi X(\xi_{j})\r)\cdot \l(\cT(\xi_j)+\cT(\xi_{j-1}) \r)+O(h^3)\\
			&=\frac{h}{2}\p_{t}\p_\xi X(\xi_{j-1})\cdot \cT(\xi_{j-1})+\frac{h}{2}\p_{t}\p_\xi X(\xi_{j})\cdot \cT(\xi_j)+O(h^3)\\
			&=\frac{h}{2}\l(-|\p_\xi X||\p_tX|^2(\xi_{j-1})-f(L)|\p_\xi X| (\xi_{j-1})\p_t X(\xi_{j-1})\cdot \cN(\xi_{j-1}) \r)\\
			&\quad +\frac{h}{2}\l(-|\p_\xi X||\p_tX|^2(\xi_{j})-f(L)|\p_\xi X| (\xi_{j})\p_t X(\xi_{j}) \cdot \cN(\xi_{j})   \r)+ O(h^3),
		\end{align*}
which together with \eqref{Qes} implies
	 \begin{equation}\label{Estimate of Remainder 2}
	 	\begin{split}
	 		\widetilde{\mathcal{R}}_j
			& = \dot{Q}_j  + \frac{Q_{j-1}+ Q_j}{4} | \dot{X}_{j-1} |^2 +\frac{Q_{j}+ Q_{j+1}}{4}| \dot{X}_j |^2\\
			&\quad +\frac{Q_{j-1}+ Q_j}{4} f(L)\dot X_{j-1}\cdot \mathcal{N}(\xi_{j-1}) +\frac{Q_{j+1}+ Q_j}{4}f(L)\dot X_{j}\cdot \mathcal{N}(\xi_{j}) = O(h^3).
	 	\end{split}
	 \end{equation}
 \item[(2).] \emph{Stability}.
 Denote  $e_j(t)=X_j(t)-x_j(t)$. Subtracting \eqref{Semidiscrete,FDM} from \eqref{Remainder term1,AP-CSF}, one gets
		\begin{align*}
			&   \dot{e}_j - \frac{2}{q_j+ q_{j+1}}\l((\cT_{j+1} - \tau_{j+1})- (\cT_j - \tau_j) \r)\\
=&  -f(L)\l(\mathcal{N}(\xi_j)-\tau_{j+1/2}^\perp \r)-\l(f(L)-f(l_h)\r)\tau_{j+1/2}^\perp\\
&\quad  + 2 \frac{  (q_j - Q_j) + (q_{j+1}-  Q_{j+1})}{(Q_j+Q_{j+1}) (q_j+q_{j+1})}\l(\cT_{j+1} - \cT_j\r)+\mathcal{R}_j
  \\
=&:  I^1_j+I^2_j + I^3_j+I^4_j.
		\end{align*}
 Multiplying both sides with $\frac{1}{2}(q_j+q_{j+1})\dot e_j$ and summing together over all $j=1,\ldots,N$, we obtain
		\[
			\frac{1}{2}\sum_{j=1}^N(q_j+q_{j+1})|\dot e_j|^2-\sum_{j=1}^N\l((\cT_{j+1} - \tau_{j+1})- (\cT_j - \tau_j) \r)\cdot \dot e_j=\sum_{k=1}^4 \sum_{j=1}^N\frac{1}{2}(q_j+q_{j+1})I_j^k\cdot \dot e_j.
		\]
Applying \eqref{FDM,induction}, Young's inequality, Assumption 2.1 and \eqref{Qes}, we arrive at
		\begin{align*}
			&-\sum_{j=1}^N\l((\cT_{j+1} - \tau_{j+1})- (\cT_j - \tau_j) \r)\cdot \dot e_j\\
			=&\ \frac{1}{2} \frac{\d }{\d t}\sum_{j=1}^Nq_j|\cT_j-\tau_{j}|^2+h\sum_{j=1}^N\l(\frac{Q_j- q_j}{Q_j} \ \delta \dot X_j \cdot ( \cT_j - \tau_j)
 +  \frac{q_j}{2Q_j} \left( \delta\dot X_j \cdot \cT_j \right)
 | \cT_j - \tau_j |^2 \r)\\
 \ge&\ \frac{1}{2}\frac{\d }{\d t}\sum_{j=1}^Nq_j|\cT_j-\tau_{j}|^2
 -C\sum_{j=1}^N\l( \frac{1}{h}(Q_j- q_j)^2 + q_j|\cT_j - \tau_j|^2
  \r),
  \end{align*}
  where for the first equality, we used the result in \cite{Deckelnick2022} (cf. page 9 in \cite{Deckelnick2022}). Employing \eqref{tausT}, \eqref{FDM,induction}, \eqref{Control of perimeter,FDM}, \eqref{Induction estimate 2} and Young's inequality, we get
  \begin{align*}
   \sum_{j=1}^N\frac{q_j+q_{j+1}}{2}I_j^1\cdot \dot e_j
  &\le \ Ch\sum_{j=1}^N\big|\mathcal{N}(\xi_j)-\tau_{j+1/2}^\perp\big||\dot e_j|
  \\
  &\le  \ C(\eps)h\sum_{j=1}^N|\cT_j-\tau_{j}|^2+\eps h\sum_{j=1}^N|\dot e_j|^2+C(\eps)h^4.
		\end{align*}
 Similarly, using \eqref{Qes}, \eqref{Xes}, Young's inequality and \eqref{Estimate of Remainder 1}, one obtains
  \begin{align*}
  \sum_{j=1}^N\frac{1}{2}(q_j+q_{j+1})I_j^3\cdot \dot e_j
  &\le \  C \sum_{j=1}^N \bigl( | Q_j - q_j| + | Q_{j+1}-q_{j+1} | \bigr) | \dot{e}_j |\\
  &\le \eps h\sum_{j=1}^N|\dot e_j|^2+\frac{C(\eps)}{h}\sum_{j=1}^N| Q_j - q_j|^2,\\
  \sum_{j=1}^N\frac{1}{2}(q_j+q_{j+1})I_j^4\cdot \dot e_j
  &\le \ \eps h\sum_{j=1}^N|\dot e_j|^2+C(\eps)h^4.
  \end{align*}
It remains to estimate the term related to $I_j^2$. Firstly we estimate the error of the perimeter by applying the trapezoidal quadrature formula and \eqref{Qes}
\begin{equation}\label{Estimate of difference of L}
	\begin{split}
		|L-l_h|
	&=\bigg|\int_{\S^1}|\p_\xi X|\ \d \xi-\sum_{j=1}^N\frac{q_j+q_{j+1}}{2} \bigg|\\
	&=\bigg|h\sum_{j=1}^N|\p_\xi X|(\xi_j)+O(h^2)-\sum_{j=1}^N\frac{q_j+q_{j+1}}{2}  \bigg|\\
	&= \bigg|\sum_{j=1}^N\frac{Q_j+Q_{j+1}}{2}+O(h^2)-\sum_{j=1}^N\frac{q_j+q_{j+1}}{2}  \bigg|\\
	&\le \sum_{j=1}^N|Q_j-q_j|+Ch^2.
	\end{split}
\end{equation}
This immediately yields
\begin{align*}
	 \sum_{j=1}^N\frac{1}{2}(q_j+q_{j+1})I_j^2\cdot \dot e_j\le &\ Ch\sum_{j=1}^N|L-l_h||\dot e_j|\le C(\eps)|L-l_h|^2+\eps h\sum_{j=1}^N|\dot e_j|^2\\
	\le &\ C(\eps)\l(\sum_{j=1}^N|Q_j-q_j|\r)^2+C(\eps)h^4+\eps h\sum_{j=1}^N|\dot e_j|^2\\
	\le &\ C(\eps)\frac{1}{h}\sum_{j=1}^N|Q_j-q_j|^2+C(\eps)h^4+\eps h\sum_{j=1}^N|\dot e_j|^2.
\end{align*}
By combining the  above inequalities, \eqref{FDM,induction} and choosing $\eps$ to be sufficiently small, we are led to
\begin{align*}
			&h\sum_{j=1}^N|\dot e_j|^2+\frac{\d }{\d t}\sum_{j=1}^Nq_j|\cT_j-\tau_{j}|^2 \le    Ch^4+C\sum_{j=1}^N\l(\frac{1}{h} (Q_j- q_j)^2 + q_j|\cT_j - \tau_j|^2 \r).
		\end{align*}
	 Through integration and utilizing Gronwall's inequality, we obtain
	  \begin{equation}\label{FDM,Stability estimate}
	 	\int^t_0h\sum_{j=1}^N|\dot e_j|^2\ \d s+\sup_{0\le s\le t} \sum_{j=1}^Nq_j|\cT_{j}-\tau_j |^2\le Ch^4+C\int^t_0\frac{1}{h} \sum_{j=1}^N (Q_j- q_j)^2\ \d s,
	 \end{equation}
	   for $0\le t\le T^*$, where $C$ is a constant depending on  $C_1, C_2, K_1(X), T$ and $f$.
\smallskip

\item[(3).] \emph{Length difference estimate}.
By using  \eqref{Induction estimate 1} and \eqref{Estimate of difference of L}, we can derive the following estimates
\begin{align*}
	&(|\dot x_j|^2-|\dot X_j|^2)\le (|\dot x_j|+|\dot X_j|)|\dot x_j-\dot X_j| \le C|\dot e_j|,\quad f(l_h)-f(L)\le  C\sum_{j=1}^N|Q_j-q_j|+Ch^2.
\end{align*}
Subtracting \eqref{Remainder term2,AP-CSF} from \eqref{FDM, length equation}, integrating from $0$ to $t$, and applying \eqref{tausT} together with the above estimate, we get
\begin{align*}
	|Q_j-q_j|(t)
	&\le \int^t_0|\dot Q_j-\dot q_j|(s)\ \d s+|Q_j- q_j|(0)\\
	&\le C\int^t_0 |q_j-Q_j|+|q_{j+1}-Q_{j+1}|+|q_{j-1}-Q_{j-1}|\ \d s\\
	&\quad + Ch\int^t_0 |\tau_j-\cT_j|+|\tau_{j+1}-\cT_{j+1}|+|\tau_{j-1}-\cT_{j-1}|\ \d s\\
	&\quad +Ch \int^t_0\sum_{j=1}^N|Q_j-q_j| \ \d s+Ch^3+Ch\int^t_0|\dot e_{j-1}|+|\dot e_{j}|\ \d s+\int^t_0|\widetilde{\mathcal{R}}_j| \ \d s.
\end{align*}
This together with \eqref{Estimate of Remainder 2} yields
\begin{align*}
	\frac{1}{h}\sum_{j=1}^N(Q_j- q_j)^2(t)
	&\le C\int^t_0h\sum_{j=1}^N|\dot e_{j}|^2\ \d s+C\int^t_0\frac{1}{h}\sum_{j=1}^N(Q_j- q_j)^2\ \d s\\
 &\qquad\qquad +C\int^t_0h\sum_{j=1}^N|\cT_j- \tau_j|^2\ \d s+Ch^4.
\end{align*}
Applying Gronwall's inequality, we get
\begin{equation}\label{FDM,Length difference estimate}
	\begin{split}
		\frac{1}{h}\sum_{j=1}^N(Q_j- q_j)^2(t)
		&\le C\int^t_0h\sum_{j=1}^N|\dot e_{j}|^2\ \d s+C\int^t_0\sum_{j=1}^Nq_j|\cT_j- \tau_j|^2\ \d s+Ch^4\\
		&\le  C\int^t_0h\sum_{j=1}^N|\dot e_{j}|^2\ \d s+C\sup_{0\le s\le t}\sum_{j=1}^Nq_j|\cT_j- \tau_j|^2+Ch^4\\
&\le Ch^4+C\int_0^t \frac{1}{h}\sum_{j=1}^N(Q_j- q_j)^2(s)ds,
	\end{split}
\end{equation}
where for the last inequality we utilized \eqref{FDM,Stability estimate}.
Hence Gronwall's inequality gives
\be\label{qQ}
\frac{1}{h}\sum_{j=1}^N(Q_j- q_j)^2(t) \le Ch^4,\quad 0\le t\le T^*.
\ee
This together with \eqref{FDM,Stability estimate} implies
\begin{equation}\label{FDM,Contiuity}
	\int^{T^*}_0 h\sum_{j=1}^N|\dot e_j|^2\ \d s+\sup_{0\le t\le T^*}\sum_{j=1}^Nq_j|\cT_j-\tau_{j}|^2 \le Ch^4.
\end{equation}
\end{enumerate}

Now we are ready to complete the proof by a continuity argument. It follows from \eqref{FDM,Contiuity} that there exists $h_2>0$ such that when $h\le h_2$,
\[| \cT_j - \tau_j |(t)
\leq h^{-\frac12} \Bigl( h \sum_{k=1}^N | \cT_k - \tau_k|^2(t)
\Bigr)^{\frac12}
\leq C h^{-\frac{1}{2}} h^2 \leq \tfrac12 h^{\frac{5}{4}},\quad 0\le t\le T^*.
\]
On the other hand, it can be easily derived from \eqref{qQ} that
 \[|Q_j(t)-q_j(t)|\le Ch^{3/2},\quad 0\le t\le T^*,\]
 which together with \eqref{Qes} yields
\[
\tfrac23 C_1 \leq q_j(t)/h \leq \tfrac32 C_2, \quad h\le h_3.
\]
By continuity we can extend $ T^*$ such that
\[ \tfrac{C_1}{2} \leq \frac{q_j(t)}{h} \leq 2 C_2, \quad \max_{j=1,\ldots,N }
| \cT_j(t)-\tau_j(t)| \leq h^{\frac54}.\]
This contradicts \eqref{FDM,induction} if $T^*<T$. Therefore, $T^*=T$. As for the estimate of $e_j$, we first notice
\[
\delta e_j=\delta X_j-\delta x_j=\frac{Q_j(\cT_j-\tau_j)}{h}+\frac{(Q_j-q_j)}{h}\tau_j.
\]
Recalling \eqref{qQ} and \eqref{FDM,Contiuity}, we immediately get
\begin{align*}
    h\sum_{j=1}^N | e_j |^2
&\le C\int^t_0 h\sum_{j=1}^N |\dot e_j |^2\d s \le Ch^4,\\
	 h\sum_{j=1}^N | \delta e_j |^2
	 &\le Ch\sum_{j=1}^N \l(|\cT_j-\tau_j|^2+\frac{(Q_j-q_j)^2}{h^2}\r)\le Ch^4,
\end{align*}
which yields
\[
	 \| X(t) - x_h(t) \|_{H^1_G}
	 =\Bigl( h \sum_{j=1}^N \bigl( | e_j |^2 + | \delta e_j |^2 \bigr) \Bigr)^{\frac{1}{2}}\le Ch^2,\quad 0\le t\le T,
\]
and the proof is completed by taking $h_0=\min\{h_1, h_2, h_3\}$.
 \qed

\section{Finite element methods}\label{sec3}
In this section, we present two finite element methods based on different formulations and establish their error estimates. The parametrization \eqref{Nonlocal,parametrization1} naturally leads to  a weak formulation: for any $v\in (H^1(\S^1))^2$, it holds
		\begin{equation}\label{Nonlocal,weak1}
	\begin{split}
			&\int_{\S^1}|\p_\xi X|\p_t X\cdot v\ \d \xi+ \int_{\S^1} \cT \cdot \p_\xi v\ \d \xi +\int_{\S^1}f(L) (\p_\xi X)^\perp\cdot v \  \d \xi=0.
		\end{split}
\end{equation}
For spatial discretization,  let $0=\xi_0<\xi_1<\ldots<\xi_N=2\pi$ be a partition of $\S^1$. We denote $h_j=\xi_j-\xi_{j-1}$ as the length of the interval $I_j:=[\xi_{j-1}, \xi_j]$ and $h=\max\limits_{j} h_j$.  We assume that the partition and the exact solution are regular in the following senses, respectively:
	  \smallskip

\textbf{(Assumption 4.1)} There exist constants $c_p$ and $c_P$ such that
	\[
	\min_{j} h_j\ge c_p h,\quad | h_{j+1}-h_j|\le c_P h^2,\quad 1\le j\le N.
	\]

	\textbf{(Assumption 4.2)} Suppose the solution of \eqref{Nonlocal,parametrization1} satisfies $X\in W^{1,\infty}\l([0,T],H^2(\S^1)\r)$, i.e., \[K_2(X):=\|X\|_{W^{1,\infty}\l([0,T],H^2(\S^1)\r)}<\infty,\]
	and there exist constants  $0<C_1<C_2$  such that \eqref{C12} holds.

  \smallskip

We define the following finite element space consisting of piecewise linear functions satisfying periodic boundary conditions:
	\[
   V_h=\l\{v\in C(\S^1,\R^2): v|_{I_j}\in P_1(I_j),\quad 1\le j\le N,\quad v(\xi_0)=v(\xi_N)\r\},
	\]
where $P_1$ denotes all polynomials with degrees at most $1$. For any continuous function $v\in C(\S^1,\R^2)$, the linear interpolation  $I_hv\in V_h$ is uniquely determined through $I_hv(\xi_j)=v(\xi_j)$ for all $1\le j\le N$ and can be explicitly written as $I_hv(\xi)=\sum\limits_{j=1}^N v(\xi_j)\varphi_j(\xi)$, where $\varphi_j$ represents the standard Lagrange basis function satisfying  $\varphi_j(\xi_i)=\delta_{ij}$.

\subsection{FEM with only the normal motion}
In this part, we  present a finite element method based on the original parametrization \eqref{Nonlocal,parametrization1}.

\begin{definition}
	We call a function
	\begin{equation}
		\label{Xh}
x_h(\xi,t)=\sum_{j=1}^Nx_j(t)\varphi_j(\xi):\S^1\times [0,T]\rightarrow \R^2
	\end{equation}
 is a semi-discrete solution of (\ref{Nonlocal,parametrization1}) if it satisfies $x_h(\xi,0)=I_hX^0$ and for all $v_h\in V_h$, it holds
		\begin{equation}\label{Semidiscrete1,weak}
		\int_{\S^1}q_h\p_t x_h\cdot v_h\ \d \xi+ \int_{\S^1} \tau_h \cdot \p_\xi v_h\ \d \xi+\int_{\S^1}\frac{\mathbf{h}^2q_h}{6}\p_\xi\p_t x_h\cdot \p_\xi v_h\mathrm{~d} \xi +\int_{\S^1}f(l_h) (\p_\xi x_h)^\perp\cdot v_h \  \d \xi=0,
		\end{equation}
 where
\be\label{qtauh}
q_h=|\p_\xi x_h|=\sum\limits_{j=1}^N \frac{q_j}{h_j}\chi_{I_j},\quad \tau_h=\frac{\p_\xi x_h}{|\p_\xi x_h|}=\sum\limits_{j=1}^N \frac{x_j-x_{j-1}}{q_j}\chi_{I_j},\ee
represent the discrete length element and unit tangent vector, respectively,
$l_h$ represents the perimeter of the evolved polygon with vertices $x_j$,
		and  $\mathbf{h}=\sum\limits_{j=1}^N h_j\chi_{I_j}$ with $\chi$ being the characteristic function.
\end{definition}

 \begin{remark}
Compared to the original formulation \eqref{Nonlocal,weak1}, here an extra term
$\int_{\S^1}\frac{\mathbf{h}^2|\p_\xi x_h|}{6}\p_\xi\p_t x_h\cdot \p_\xi v_h\mathrm{~d} \xi$ is introduced in \eqref{Semidiscrete1,weak}, which reduces to the so-called mass-lumped scheme \eqref{Semidiscrete1,lump mass}. Clearly this
term does not affect the convergence order for a linear finite element method.
As was interpreted in \cite{Jiang2022,Dziuk1999}, this mass-lumped version can preserve the length shortening property for the CSF/AP-CSF, which was missing for the original formula.
 \end{remark}

Taking $v_h=(\varphi_j, 0)$ and $v_h=(0, \varphi_j)$ for $j=1,\ldots, N$ in  \eqref{Semidiscrete1,weak}, we are led to the following $2N$ ordinary differential equations:
\begin{equation}\label{Semidiscrete1,lump mass}
	\frac{q_j+q_{j+1}}{2}\dot x_j=\tau_{j+1}-\tau_{j}-f(l_h)(x_{j+1}-x_{j-1})^\perp,
\end{equation}
where $\tau_{j}$ is the discrete tangent defined as \eqref{qtau}. Furthermore, we have the following identities
\begin{align}
	\dot q_j		&=-\frac{1}{q_j+q_{j+1}}|\tau_{j+1}-\tau_{j}|^2-\frac{1}{q_j+q_{j-1}}|\tau_{j-1}
-\tau_{j}|^2 +\tau_j\cdot\l(r_j-r_{j-1}  \r) \label{dq1}\\
		&=-\frac{q_j+q_{j+1}}{4}|\dot x_j-r_j|^2-\frac{q_j+q_{j-1}}{4}|\dot x_{j-1}-r_{j-1}|^2+\tau_j\cdot\l(r_j-r_{j-1}  \r),\label{dq2}
			\end{align}
	where for simplicity we denote
	\begin{equation}\label{rrj}
	r_j=-f(l_h)\frac{n_{j}q_j+n_{j+1}q_{j+1}}{q_j+q_{j+1}}.
	\end{equation}

\begin{theorem}\label{Error estimate 1}
Let $X(\xi,t)$ be a solution of \eqref{Nonlocal,parametrization1}  satisfying Assumption 4.2. Assume that the partition of $\S^{1}$ satisfies Assumption 4.1. Then there exists $h_0>0$ such that for all $0<h\le h_0$, there exists a unique semi-discrete solution $x_{h}$ for \eqref{Semidiscrete1,weak}. Furthermore, the solution satisfies
	\begin{equation}\label{Main estimate 1}
	\begin{split}
		\int^T_0\|\p_t X-\p_t x_h\|^2_{L^2}\mathrm{d}t+\sup_{t\in [0,T]}\|X-x_h\|^2_{H^1}\le Ch^2,
	\end{split}
	\end{equation}
where $h_0$ and $C$ depend on $c_p, c_P,C_1,C_2, T$, $K_2(X)$ and $f$.
	
\end{theorem}

Before presenting the proof of Theorem \ref{Error estimate 1}, we first list a lemma which will be used later.

\begin{lemma}\cite[Lemma 4.2]{Jiang2022}\label{Control lemma}
Under Assumptions 3.1 and 3.2, suppose further
   \[
    \int_{\S^{1}}\left|\cT-\tau_h\right|^{2}q_h\d \xi+\||\p_\xi X|-q_h\|_{L^2}^{2}\leq Ch^{2},\quad \forall\ t\in[0,T^*],
   \]
then  there exists a constant $h_0$ such that for any $0<h\le h_0$, we have
   \[
  \inf_{\xi}q_h\ge 3C_1/4,\quad \text{and}\quad \sup_{\xi}q_h\le 3C_2/2,\quad \forall\ t\in[0,T^*],
   \]
   where $C_1$ and $C_2$ are the lower and upper bounds of $|\p_\xi X|$ shown in \eqref{C12}.
\end{lemma}

\Proof Similar to the proof of Theorem \ref{Error estimate 3}, we apply the continuity argument. Define
\begin{equation}\label{Continuous method}
		T^*=\sup\{t\in [0,T]: \eqref{Semidiscrete1,weak}\,\,\, \text{has a unique solution}\,\,\,x_h \,\,\,\text{and}\,\,\,\inf q_h\ge C_1/2,\,\,\sup q_h\le 2C_2 \}.
\end{equation}
Since the nonlinear terms in \eqref{Semidiscrete1,lump mass} are locally Lipschitz with respect to $x_j$, the local existence and uniqueness follow from  standard ODE theory, and thus $T^*>0$. Moreover, due to the Lipschitz property of $f$ and Assumption \eqref{Continuous method}, for any  $ t\in [0,T^*]$, it holds  that
 \begin{equation}\label{Control of perimeter}
 \begin{split}
 	2\pi C_1\le  L\le 2\pi C_2 ,&\quad \pi C_1 \le l_h\le 4\pi C_2 ,\quad  |f(L)|\le C,
 	 \end{split}
 \end{equation}
 where $C$ is a constant depending on $C_1, C_2, f$.
\begin{enumerate}
	\item[(1)] \emph{Stability}.
Taking the difference between  \eqref{Nonlocal,weak1} and \eqref{Semidiscrete1,weak}, and choosing   $v_h=I_h(\p_tX)-\p_t x_h\in V_h$, we get
	  \begin{align*}
  	&\quad\int_{\S^1}|\p_t X-\p_t x_{h} |^2q_h \d\xi +\int_{\S^1}\l(\cT-\tau_h  \r)\l(\p_\xi\p_t X- \p_\xi\p_t x_h\r) \d\xi\\
  	&=\int_{\S^1}\p_t X\cdot (q_h-|\p_\xi X|)\l(I_h\p_t X-\p_t x_{h} \r)\d\xi+\int_{\S^1}\frac{\mathbf{h}^2q_h}{6}\p_\xi\p_t x_h\cdot \p_\xi \l(I_h\p_t X-\p_t x_{h} \r)\d\xi \\
  	&\quad +\int_{\S^1}q_h(\p_t X-\p_t x_h)\cdot(\p_t X-I_h\p_t X) d\xi+\int_{\S^1}\l(\cT-\tau_h  \r)\cdot \l(\p_\xi\p_t X- \p_\xi I_h\p_t X\r)\d \xi \\
  	&\quad +\int_{\S^1}f(L)\l(\p_\xi X-\p_\xi x_h \r)^\perp\cdot \l(\p_t x_h-I_h\p_t X\r)\d\xi\\
  	&\quad + \int_{\S^1}\Big(f(L)-f(l_h) \Big)\l(\p_\xi x_h \r)^\perp \cdot \l(\p_t x_h-I_h\p_t X\r)\d \xi 	=:\ J_1+J_2+J_3+J_4+J_5+J_6.
  \end{align*}
    The estimates of the second term on the left side and $J_j$ for $1\le j\le 4$ can be found in \cite[Lemma 5.1]{Dziuk1999} or \cite[Lemma 4.1]{Jiang2022}, which can be summarized as follows:
  	\begin{align*}
  &\quad\int_{\S^1}\l(\cT-\tau_h \r)\cdot\l(\p_\xi\p_t X- \p_\xi\p_t x_h\r)  \d\xi\\
&\ge \frac{1}{2}\frac{\d}{\d t}\Big(\int_{\S^1}|\mathcal{T}-\tau_h|^2q_h  \d\xi\Big)-C\|\p_\xi\p_tX\|_{L^\infty}\Big(\int_{\S^1}|\cT-\tau_h|^2 q_h d\xi
+\||\p_\xi X|-q_h\|_{L^2}^2\Big),\\
  &J_1+J_2+J_3+J_4\le \varepsilon\int_{\mathbb{S}^1}\left|\partial_t X-\partial_t X_{h}\right|^{2}q_h \mathrm{d} \xi+C(\varepsilon)\left\|\partial_t X\right\|_{L^{\infty}}^{2} \||\p_\xi X|-q_h\|_{L^2}^{2} \\
  &\qquad\qquad\qquad\qquad\qquad\qquad +C(\varepsilon)h^2\|\partial_tX\|_{H^1}^2+C\int_{\mathbb{S}^1}|\mathcal{T}-\tau_h|^2 q_h\mathrm{d}\xi,
  	\end{align*}
where $\varepsilon$ is a generic positive constant which will be chosen later.
 For $J_5$ and $J_6$, in view of the Lipschitz property of $f$, \eqref{Control of perimeter}, and the identity
  \be\label{xxh}
  |\partial_\xi X-\partial_\xi x_h |^2=
               (|\p_\xi X|-q_h)^2+|\p_\xi X|q_h |\mathcal{T}-\tau_h|^2,
               \ee
  applying similar techniques in \cite{Jiang2022} (cf. proof of Lemma 4.1), we can get
  \begin{align*}
    J_5
    &=\int_{\S^1}f(L)\l(\p_\xi X-\p_\xi x_h \r)^\perp\cdot \l(\p_t X-I_h\p_t X\r) \d\xi\\
    &\qquad\qquad\qquad +\int_{\S^1}f(L)\l(\p_\xi X-\p_\xi x_h \r)^\perp\cdot \l(\p_t x_h-\p_t X\r)  \d \xi\\
    &\le C\l\|\p_\xi X\r\|_{L^\infty}\int_{\S^1}\l|\cT-\tau_h \r|^2q_h \d \xi + C\int_{\S^1}(|\p_\xi X|-q_h)^2 \d \xi+Ch^2\|\p_t X\|_{H^1}^2\\
    &\ +C(\varepsilon )\l\|\p_\xi X\r\|_{L^\infty}\int_{\S^1}\l|\cT-\tau_h \r|^2q_h \d \xi + C(\varepsilon )\||\p_\xi X|-q_h\|_{L^2}^2+\varepsilon  \int_{\S^1}\l|\p_t x_h-\p_t X \r|^2q_h\d \xi,\\
  	J_6
  	&=\int_{\S^1}\big(f(L)-f(l_h) \big)\l(\p_\xi x_h \r)^\perp \cdot \l(\p_t X-I_h\p_t X\r)\d \xi\\
   &\qquad\qquad\qquad +\int_{\S^1}\big(f(L)-f(l_h)  \big)\l(\p_\xi x_h \r)^\perp \cdot \l(\p_t x_h-\p_t X\r)\d \xi\\
  	&\le C|L-l_h|^2+ C\|\p_t X-I_h\p_t X\|_{L^2}^2 +C(\varepsilon )|L-l_h|^2+\varepsilon \int_{\S^1}q_h\l|\p_t x_h-\p_t X \r|^2\ \d \xi\\
  	&\le C(\varepsilon)\||\p_\xi X|-q_h\|^2_{L^2}+Ch^2\|\p_t X\|^2_{H^1}+\varepsilon \int_{\S^1}q_h\l|\p_t x_h-\p_t X \r|^2\ \d \xi.
  \end{align*}	
Here we use  the inequalities $|L-l_h|^2\le \||\p_\xi X|-q_h\|_{L^1}^2 \le C\||\p_\xi X|-q_h\|_{L^2}^2$. Combining all the above estimates, we are led to
  \begin{align*}
  	&\int_{\S^1}|\p_t X-\p_t x_{h} |^2q_h \d\xi +\frac{1}{2}\frac{\d}{\d t}\int_{\S^1}|\mathcal{T}-\tau_h|^2q_h \d\xi\le 4\eps \int_{\S^1}|\p_t X-\p_t x_{h} |^2q_h \d\xi\\
  	&\qquad+C(\eps)h^2\|\p_t X\|^2_{H^2} +C(\eps,K_2(X)) \||\p_\xi X|-q_h\|_{L^2}^2+C(\eps,K_2(X)) \int_{\S^1}|\mathcal{T}-\tau_h|^2q_h\d\xi.
  \end{align*}
 Choosing $\eps$ small enough, integrating both sides with respect to time from $0$ to $t$ and applying Gronwall's argument, we arrive at
  \begin{equation}\label{Stability estimate}
  \int_{0}^{t} \int_{\S^{1}}|\p_t X-\p_t x_h|^{2}q_h \d \xi \d s +\sup\limits_{0\le s\le t}\int_{\S^{1}}|\cT-\tau_h|^{2}q_h\d \xi\leq C \int_{0}^{t} \||\p_\xi X|-q_h\|_{L^2}^{2}\d s +C h^{2},
    \end{equation}
    where $C$ is a constant depending on $c_p$, $c_P$, $C_1$, $C_2$, $T$, $K_2(X)$ and $f$.
    \item[(2)]  \emph{Length difference estimate}.
     Applying the Lipschitz property of $f$ and \eqref{Control of perimeter}, a mild modification of the proof of
        \cite[Lemma 4.3, Lemma 4.4]{Jiang2022} enables us to establish the same length difference estimate as in \cite{Jiang2022}:
\begin{equation}\label{Length difference estimate}
		\||\p_\xi X|-q_h\|_{L^2}^2
	\le C\int_{0}^{t} \int_{\S^{1}}\left|\p_t X-\p_t x_h\right|^{2}q_h \d \xi \d s+C\int^t_0\int_{\S^{1}}\left|\cT-\tau_h\right|^{2}q_h\d \xi\d s+Ch^{2},
\end{equation}
  where $C$ depends on $c_p$, $c_P$, $C_1$, $C_2$, $T$, $K_2(X)$ and $f$. For the details, we refer to \cite{Jiang2022}.
\end{enumerate}

Combining \eqref{Stability estimate} and \eqref{Length difference estimate}, employing Gronwall's inequality, we derive
\be\label{temp1}
	\int_{0}^{t} \int_{\S^{1}}\left|\p_t X-\p_t x_h\right|^{2}q_h \d \xi \d s+\sup_{0\le s\le t} \int_{\S^{1}}\left|\cT-\tau_h\right|^{2}q_h\d \xi \le Ch^2,\quad \forall\ t\in [0,T^*],
\ee
which together with \eqref{Length difference estimate} yields
\[
	 \||\p_\xi X|-q_h\|_{L^2}^2\leq Ch^{2}.
\]
Applying Lemma \ref{Control lemma}, there exists $h_0>0$ depending  on  $c_p, c_P, C_1, C_2, T$, $K_2(X)$ such that for any $0<h\le h_0$, we have
  \[
\inf q_h\ge 3C_1/4,\quad \text{and}\quad \sup q_h\le 3C_2/2,\quad  t\in [0, T^*].
  \]
By standard ODE theory, we can uniquely  extend the above semi-discrete solution in a neighborhood of $T^*$, and thus $T^*=T$.
The estimate \eqref{Main estimate 1} can be concluded similarly as in \cite[Theorem 2.5]{Jiang2022} by integration, \eqref{xxh} and \eqref{temp1}:
\begin{align*}
&\quad\|X(\cdot, t)-x_h(\cdot, t)\|_{H^1}^2=
\int_{\S^1}|X-x_h|^2\d\xi+\int_{\S^1}|\p_\xi X-\p_\xi x_h|^2\d\xi\\
&\le 2\int_{\S^1}\big(\int_0^t \p_t X-\p_t x_h\d s\big)^2\d\xi+2\|X^0-I_h X^0\|_{L^2}^2+\||\p_\xi X|-q_h\|_{L^2}^2+\int_{\S^1}|\mathcal{T}-\tau_h|^2|\p_\xi X|q_h \d\xi\\
&\le 2\int_{\S^1} T\int_0^t |\p_t X-\p_t x_h|^2\d s\d\xi+Ch^2
\le Ch^2,
\end{align*}
and the proof is completed.
\qed

\subsection{FEM with tangential motions}\label{sec4}
The aforementioned methods are developed based on the equation \eqref{Nonlocal} and only normal motion is allowed. They might suffer from the fact that the mesh will have inhomogeneous properties during the evolution, for instance, some nodes may cluster and the mesh may become distorted. This will lead to instability and even the breakdown of the simulation. To address this challenge, various techniques have been proposed to improve the mesh quality for evolving various types of geometric flows in the literature, such as mesh redistribution \cite{Bansch2005}, and the introduction of artificial tangential velocity \cite{Barrett2020,Li2024,Mikula2004b,Sevcovic2001}.

In this subsection, to achieve equipartition property for long-time evolution, we derive another formulation of \eqref{Nonlocal} by introducing a tangential velocity.  We consider the equation
\begin{align*}
	\p_t X = (\kappa-f(L))\cN + \gamma(X) \cT,
\end{align*}
where $\cN,\cT$ are the unit normal vector and tangent vector respectively, and $\gamma$ is the tangential velocity to be determined. It is important to note that the presence of tangential velocity has no impact on the shape of evolving curves \cite{Deckelnick1995,Elliott2017}, and suitable choices of tangential velocity may help the redistribution of mesh points \cite{Mikula2004,Mikula2004b,Sevcovic2001}. As mentioned in the introduction, inspired by the work of \cite{Deckelnick1995,Elliott2017} for curve shortening flow, we consider an explicit tangential velocity given by
\[
	\gamma(X)=\frac{\p_\xi X\cdot \p_{\xi\xi} X }{|\p_\xi X|^3}.
\]
More generally, for a fixed parameter $0< \alpha\le 1$, we consider a series of reparametrizations $X_\alpha$ which are determined by
\begin{equation}\label{Nonlocal, tangent motion}
	\alpha \p_t X_\alpha +(1-\alpha) (\p_t X_\alpha\cdot \cN)\cN
	= \frac{\p_{\xi\xi}X_\alpha}{|\p_\xi X_\alpha |^2}  -f(L)\cN;\quad X_\alpha(\xi,0)=X^0(\xi).
\end{equation}
Below we provide three justifications for \eqref{Nonlocal, tangent motion}.
\begin{itemize}
	\item[(i)] The solution $X_\alpha$ of  the evolution equation   \eqref{Nonlocal, tangent motion} has the same shape as the standard parametrization equation  \eqref{Nonlocal,parametrization1} since they share the same normal velocity
	\begin{align*}
		\p_{t}X_\alpha\cdot \cN
		&= \alpha \p_{t}X_\alpha\cdot \cN+(1-\alpha)\p_{t}X_\alpha\cdot \cN\\
		&=\l(\frac{\p_{\xi\xi}X_\alpha}{|\p_\xi X_\alpha |^2}-\gamma(X_\alpha)\cT  -f(L)\cN\r)\cdot \cN\\
		&=\l(\frac{1}{|\p_\xi X_\alpha |} \p_\xi\l(\frac{\p_\xi X_\alpha}{|\p_\xi X_\alpha|} \r) -f(L)\cN\r)\cdot \cN\\
		&=\kappa-f(L),
	\end{align*}
where we note that the curvature $\kappa$ and perimeter $L$ are geometric quantities independent of the parametrization.

	\item[(ii)] The evolution of \eqref{Nonlocal, tangent motion} has asymptotic equidistribution property in a continuous level. More precisely, suppose $X_\alpha^e$ is the equilibrium of \eqref{Nonlocal, tangent motion}, i.e., $\p_tX_\alpha^e=0$, then formally we have
	\begin{align*}
		\p_\xi|\p_\xi X_\alpha^e|
		&=\p_{\xi\xi} X_\alpha^e\cdot \cT=\l(\frac{\p_{\xi\xi} X_\alpha^e}{|\p_{\xi} X_\alpha^e|^2}-f(L)\cN\r) \cdot |\p_{\xi} X_\alpha^e|^2\cT=0,
	\end{align*}
	which means the equilibrium has constant arc-length. This leads us to expect that the corresponding numerical solution for  \eqref{Nonlocal, tangent motion} has equidistributed mesh points for long-time evolution.

	\item[(iii)] As explained in \cite[Section 8]{Elliott2017}, we can write  the standard parametrization equation  \eqref{Nonlocal,parametrization1} as
	\[
		\p_tX=\Delta_{\Gamma[X]}X-f(L)\cN,
	\]
	where $\Gamma[X]$ is the image of $X$ and $\Delta_{\Gamma[X]}$ is the Laplace-Beltrami operator over the curve $\Gamma[X]$. The DeTurck's trick for operator $\Delta_{\Gamma[X]}$ maintains  the normal term $f(L)\cN$ unaffected and leads to \eqref{Nonlocal, tangent motion}. In this aspect, the nonlocal flows can be viewed as a natural generalization of \cite[Section 8]{Elliott2017}.
\end{itemize}

Next we present a finite element method for \eqref{Nonlocal, tangent motion}. For fixed $\alpha$, multiplying  $|\p_\xi X |^2$ for both  sides of \eqref{Nonlocal, tangent motion} (below we omit the subscript $\alpha$ for simplicity), we obtain the following weak formulation: for any $\ v\in (H^1(\S^1))^2$, it holds
\begin{equation}\label{Nonlocal,weak2}
\int_{\S^1}|\p_\xi X |^2(\alpha \p_t X +(1-\alpha) (\p_t X\cdot \cN)\cN)\cdot  v \d \xi +\int_{\S^1}\p_\xi X\cdot \p_\xi v \d \xi +\int_{\S^1} f(L)|\p_\xi X|^2\cN\cdot  v \d \xi =0.
\end{equation}
We use the same spatial discretization for $\S^1$ as in the last subsection and assume it satisfies Assumption 4.1. We further assume the  exact solution of \eqref{Nonlocal, tangent motion} is regular in the following sense.

\smallskip

\textbf{(Assumption 4.3)} Suppose that  the  solution of \eqref{Nonlocal, tangent motion} with an initial value $X^0\in H^{2}(\S^1)$ satisfies $X\in W^{1,\infty}\l([0,T],H^2(\S^1)\r)$, i.e., \[K_2(X):=\|X\|_{W^{1,\infty}\l([0,T],H^2(\S^1)\r)}<\infty,\]
and there exist constants  $0<C_1<C_2$  such that \eqref{C12} holds.

\begin{definition}
	We call a function $	x_h(\xi,t)=\sum\limits_{j=1}^Nx_j(t)\varphi_j(\xi):\S^1\times [0,T]\rightarrow \R^2$
 is a semi-discrete solution of (\ref{Nonlocal, tangent motion}) if it satisfies $x_h(\xi,0)=I_hX^0$  and
			\begin{equation}\label{Semidiscrete2,weak}
\int_{\S^1}|\p_\xi x_h |^2(\alpha \p_t x_h +(1-\alpha) (\p_t x_h\cdot n_h)n_h)\cdot  v_h \ \d \xi +\int_{\S^1}\p_\xi x_h\cdot \p_\xi v_h \ \d \xi \\
	 +\int_{\S^1} f(l_h)|\p_\xi x_h|^2n_h\cdot  v_h \ \d \xi =0,
\end{equation}
for any $v_h\in V_h$, where $n_h=\tau_h^\perp$ represents the piecewise  unit normal vector.
\end{definition}

  \smallskip

\begin{theorem}\label{Error estimate 2}
Let $X(\xi,t)$ be a solution of \eqref{Nonlocal, tangent motion}  satisfying Assumption 4.3. Assume that the partition of $\S^{1}$ satisfies Assumption 4.1. Then there exists $h_0>0$ such that for all $0<h\le h_0$, there exists a unique semi-discrete solution $x_{h}$ for \eqref{Semidiscrete2,weak}. Furthermore, the solution satisfies
\begin{align*}
    \sup_{t\in[0,T]}|X-x_h |_{H^1}^2 +\alpha \int^T_0\|\p_t X-\p_t x_h \|_{L^2}^2\ \d t&+(1-\alpha ) \int^T_0\|n_h\cdot(\p_t X-\p_t x_h) \|_{L^2}^2\ \d t \\
    &\qquad\qquad\qquad \le {CTh^2+Ce^{MT/\alpha}h^2},
\end{align*}
where $h_0$, $C$ and $M$ depend on $c_p$, $c_P$, $C_1$, $C_2$, $K_2(X)$, and $f$.

\end{theorem}

\Proof
	Fix $\alpha\in(0,1]$. We consider a Banach space $Z_h=C([0,T],V_h)$ equipped with the norm
\[
\|v_h\|_{Z_h}:=\sup_{t\in[0,T]} \|v_h(t)\|_{L^2},\quad v_h\in Z_h,
\]
and a nonempty closed  convex subset $B_h$  of $Z_h$ defined by
\begin{equation}\label{Definition of B_h}
	B_h:=\bigg\{v_h\in Z_h|\sup_{t\in[0,T]}e^{-M t/\alpha}\|(\p_\xi X-\p_\xi v_h )(t) \|_{L^2}^2\le K^2h^2\ \text{and}\ v_h(\cdot,0)=I_hX^0(\cdot) \bigg\},
\end{equation}
where $M,K>0$ are constants that will be determined later. For any $u_h\in B_h$, applying  interpolation error, inverse inequality and \eqref{Definition of B_h}, one can easily derive
	\begin{align*}
		\|(\p_\xi X-\p_\xi u_h)(t) \|_{L^\infty}
		&\le \|(\p_\xi X- I_h\p_\xi X)(t) \|_{L^\infty}+ \|( I_h\p_\xi X-\p_\xi I_hu_h)(t) \|_{L^\infty}\\
		&\le Ch^{1/2}+ Ch^{-1/2}\|( I_h\p_\xi X-\p_\xi u_h)(t) \|_{L^2}\\
		&\le Ch^{1/2}+ Ch^{-1/2}\l(\|(\p_\xi X- I_h\p_\xi X)(t) \|_{L^2}+\|(\p_\xi X-\p_\xi u_h)(t) \|_{L^2} \r)\\
		&\le Ch^{1/2}\l(1+e^{\frac{Mt}{2\alpha}}K\r).
	\end{align*}
	It follows from Assumption 4.3 that there  exists a constant $h_0>0$ depending on $\alpha,M,$\ $K,T,K_2(X)$  such that for any $0<h\le h_0$, we have
	\begin{equation}\label{Control of length element 2}
		\inf_{\xi}|\p_\xi u_h|\ge C_1/2,\quad \sup_{\xi}|\p_\xi u_h |\le 2C_2 .
	\end{equation}
 Setting $\widehat{q_h}=|\p_\xi u_h|$ and denoting $\widehat{l_h}$ as the perimeter of $u_h$, due to the Lipschitz property of $f$, it holds that
 \begin{equation}\label{Control of perimeter 2}
 \begin{split}
 	2\pi C_1\le  L\le 2\pi C_2 ,&\quad \pi C_1 \le \widehat{l_h}\le 4\pi C_2 ,\quad |f(\widehat{l_h})|\le C,
 \end{split}
 \end{equation}
 where $C$ is a constant depending  on $C_1, C_2$ and $f$.  We define a continuous  map  $F:B_h\rightarrow Z_h$ as follows. For  any $u_h\in B_h$, we define $y_h$ as the unique solution of the following linear equation:
	\begin{equation}\label{Weak, ODE}
	\begin{split}
		&\int_{\S^1}\widehat{q_h}^2 (\alpha\p_t y_h+(1-\alpha)(\p_ty_h\cdot \widehat{n_h})\widehat{n_h} )\cdot v_h \ \d \xi +\int_{\S^1}\p_\xi y_h\cdot \p_\xi v_h \ \d \xi \\
		&\qquad \qquad +\int_{\S^1} f(\widehat{l_h})\widehat{q_h}(\p_\xi y_h)^\perp \cdot  v_h \ \d \xi =0,\quad \forall \ v_h\in V_h,
	\end{split}		
	\end{equation}
	with initial data $y_h( 0)=I_hX^0$, where  $\widehat{n_h}=\l(\frac{\p_\xi u_h}{|\p_\xi u_h|}\r)^{\perp}$.
\begin{enumerate}
	\item[(1)]  \emph{Length difference estimate for $u_h\in B_h$}.  Applying \eqref{Definition of B_h} and the triangle inequality, we obtain
	\begin{equation}\label{Basic estimate, tangent motion}
		\|(|\p_\xi X|-\widehat{q_h})(t) \|_{L^2}\le \|(\p_\xi X-\p_\xi u_h)(t)\|_{L^2}\le Khe^{Mt/(2\alpha)},\quad 0\le t\le T.
	\end{equation}
	\item[(2)] \emph{Stability estimate for $y_h\in Z_h$}. Taking $v=v_h$ in \eqref{Nonlocal,weak2} and  subtracting \eqref{Weak, ODE} from \eqref{Nonlocal,weak2}, we get
	\begin{align*}
		&\int_{\S^1}\widehat{q_h}^2\l(\alpha \l(\p_t X-\p_t y_h\r)\cdot  v_h+(1-\alpha)\l(\p_t X-\p_t y_h\r)\cdot \widehat{n_h} (\widehat{n_h}\cdot v_h) \r) \ \d \xi\\
  &\qquad\qquad\qquad + \int_{\S^1}\l(\p_\xi X- \p_\xi y_h\r)\cdot \p_\xi v_h \ \d \xi\\
		=&\ \int_{\S^1}(\widehat{q_h}^2-|\p_\xi X|^2 )(\alpha\p_t X\cdot v_h+(1-\alpha)(\p_t X\cdot \widehat{n_h})(\widehat{n_h}\cdot v_h) )\ \d \xi  \\
		& +(1-\alpha)\int_{\S^1}|\p_\xi X|^2(\p_t X\cdot (\widehat{n_h}-\cN)(\widehat{n_h}\cdot v_h)+(\p_t X\cdot \cN)(\widehat{n_h}-\cN)\cdot v_h )  \ \d\xi \\
		&+\int_{\S^1}\l(-f(L)|\p_\xi X|+f(\widehat{l_h})\widehat{q_h}\r)(\p_\xi X)^\perp\cdot v_h \ \d \xi\\
		& -\int_{\S^1} f(\widehat{l_h})\widehat{q_h}\l((\p_\xi X)^\perp-(\p_\xi y_h)^\perp \r)\cdot v_h \ \d \xi =:\ J_1+J_2+J_3+J_4.
	\end{align*}
Choosing  $v_h=I_h(\p_tX)-\p_ty_h\in V_h$, the estimates of the left-hand hand side and $J_1,J_2$ can be found in  \cite[(3.7)]{Elliott2017}, which can be summarized as
	\begin{equation}\label{Elliott-Fritz estimate}
	\begin{split}
		&\quad\frac{\d }{\d t}\|\p_\xi X-\p_\xi y_h\|^2_{L^2}+\frac{C_1^2}{8}\alpha\|\p_t X-\p_t y_h \|^2_{L^2}+\frac{C_1^2}{4}(1-\alpha)\|\widehat{n_h}\cdot (\p_t X-\p_t y_h) \|^2_{L^2}\\
    	&\le \ Ch^2(1+\|\p_t X\|^2_{H^{2}})+\|\p_\xi X-\p_\xi y_h\|^2_{L^2}+Ch^2K^2 e^{Mt/\alpha}/\alpha+2|J_3|+2|J_4|.
	\end{split}		
	\end{equation}
For the terms of $J_3$ and $J_4$, in view of the Lipschitz property of $f$ and the inequality
\[|L-\widehat{l_h} |\le C\||\p_\xi X|-\widehat{q_h} \|_{L^2},\]
applying \eqref{Control of length element 2}, \eqref{Control of perimeter 2},  \eqref{Basic estimate, tangent motion} and Young's inequality, we get
    \begin{align*}
	|J_3|
	&=\bigg|\int_{\S^1}\l(-f(L)+f(\widehat{l_h}) \r)|\p_\xi X| (\p_\xi X)^\perp\cdot v_h \ \d \xi -\int_{\S^1}f(\widehat{l_h}) \l(|\p_\xi X|-\widehat{q_h} \r) (\p_\xi X)^\perp\cdot v_h \ \d \xi\bigg|\\
	&\le C\int_{\S^1}|L-\widehat{l_h}|\l( |I_h(\p_tX)-\p_t X |+|\p_t X-\p_t y_h |\r)\ \d\xi\\
	&\qquad +C\int_{\S^1}||\p_\xi X|-\widehat{q_h}|  \l( |I_h(\p_tX)-\p_t X |+|\p_t X-\p_t y_h |\r)\ \d\xi\\
	&\le Ch|L-\widehat{l_h}|\|\p_t X\|_{H^1} +C|L-\widehat{l_h}|\|\p_t X-\p_t y_h \|_{L^2} \\
	&\qquad +Ch\||\p_\xi X|-\widehat{q_h} \|_{L^2}\|\p_t X\|_{H^1}+C\||\p_\xi X|-\widehat{q_h}\|_{L^2}  \|\p_t X-\p_t y_h \|_{L^2}\\
	&\le CKe^{\frac{Mt}{2\alpha}}h\|\p_t X-\p_t y_h \|_{L^2}+CKe^{\frac{Mt}{2\alpha}}h^2\|\p_t X\|_{H^1}\\
	&\le  \frac{C(\varepsilon) e^{Mt/\alpha}K^2h^2}{4\alpha} +\varepsilon \alpha \|\p_t X-\p_t y_h \|_{L^2}^2+\alpha h^2\|\p_t X\|_{H^1}^2,
\end{align*}
and
\begin{align*}
		|J_4|
	&=\bigg|\int_{\S^1} f(\widehat{l_h})\widehat{q_h}\l((\p_\xi X)^\perp-(\p_\xi y_h)^\perp \r)\cdot v_h\ \d \xi\bigg|\\
	&\le Ch\| \p_\xi X-\p_\xi y_h \|_{L^2}\|\p_t X\|_{H^1}+C\| \p_\xi X-\p_\xi y_h \|_{L^2} \|\p_t X-\p_t y_h \|_{L^2} \\
	&\le   \| \p_\xi X-\p_\xi y_h \|_{L^2}^2+C h^2 \|\p_t X\|_{H^1}^2 +\frac{C(\varepsilon )}{4 \alpha} \| \p_\xi X-\p_\xi y_h \|_{L^2}^2  +  \varepsilon\alpha  \|\p_t X-\p_t y_h \|_{L^2}^2.
\end{align*}
Combining all the above estimate and taking $\varepsilon$ small enough,  we obtain
\begin{equation}\label{Combination, tangent motion}
\begin{split}
	&\quad\frac{\d }{\d t}\|\p_\xi X- \p_\xi y_h \|_{L^2}^2+\frac{C_1^2}{16}\alpha \|\p_t X-\p_t y_h\|_{L^2}^2+\frac{C_1^2}{4}(1-\alpha)\|\widehat{n_h}\cdot(\p_t X-\p_ty_h ) \|_{L^2}^2  \\
	&\le \ Ch^2+C(1+1/\alpha)\|\p_\xi X-\p_\xi y_h\|^2_{L^2}+Ch^2K^2 e^{Mt/\alpha}/\alpha,
\end{split}	
\end{equation}
where $C$ depends on $c_p$, $c_P$, $C_1$, $C_2$, $T$, $K_2(X)$, and $f$.
This directly gives
\[\frac{\d }{\d t}\|\p_\xi X- \p_\xi y_h \|_{L^2}^2\le \ Ch^2+C(1+1/\alpha)\|\p_\xi X-\p_\xi y_h\|^2_{L^2}+Ch^2K^2 e^{Mt/\alpha}/\alpha.\]
Thus we get
\begin{align*}
&\|\p_\xi X(t)- \p_\xi y_h(t)\|_{L^2}^2\\
&\le \|\p_\xi X^0- \p_\xi y_h(0)\|_{L^2}^2 e^{C(1+\frac1{\alpha})t}+Ch^2\int_0^t e^{C(1+\frac1{\alpha})(t-s)}
\l(1+K^2e^{Ms/\alpha}/\alpha\r)ds\\
&\le Ce^{C(1+\frac1{\alpha})t}h^2+CK^2h^2\frac{e^{\frac{M}{\alpha}t}
-e^{C(1+\frac{1}{\alpha})t}}{M-C(1+\alpha)},
\end{align*}
which yields
\begin{equation}\label{Induction estimate, tangent motion}
	e^{-Mt/\alpha}\|\p_\xi X-\p_\xi y_h\|^2_{L^2}
\le Ch^2 e^{(-\frac{M}{\alpha}+\frac{C}{\alpha}+C)t}+\frac{CK^2h^2}{M-C(1+\alpha)}
\le K^2h^2,
\end{equation}
if we select $M\ge 3C+C\alpha$ and $K^2\ge 2C$. Hence, by plugging  \eqref{Induction estimate, tangent motion} into \eqref{Combination, tangent motion}, integrating  from $0$ to $T$,  we arrive at
\begin{align}
   	\max_{t\in [0,T]}&\|\p_\xi X-\p_\xi y_h\|^2_{L^2}+\alpha\int^T_0\|\p_t X-\p_t y_h\|_{L^2}^2\ \d t+(1-\alpha)\int^T_0\|\widehat{n_h}\cdot(\p_t X-\p_ty_h ) \|_{L^2}^2\d t \notag\\
   &\le C(1+T)h^2+CK^2h^2(1+1/\alpha)\int_0^T e^{Ms/\alpha}ds\notag\\
   &\le C(1+T)h^2+K^2h^2(e^{MT/\alpha}-1)
   \le CTh^2+Ce^{MT/\alpha}h^2,\label{Stability estimate for y_h}
\end{align}
where $C$ and $M$ are constants depending on $c_p,c_P,K_2(X)$, $C_1,C_2$ and $f$.
\end{enumerate}

Now we complete the proof by applying  Schauder's fixed point theorem for $F$. Indeed, it follows from assumption \eqref{Definition of B_h} and \eqref{Induction estimate, tangent motion} that $F(B_h)\subset B_h$. Furthermore, it can be easily derived from \eqref{Stability estimate for y_h} and the assumption $y_h(0)=I_hX^0$ that $\|y_h\|_{W^{1,2}([0,T],V_h)}\le C$, which, together with the Sobolev embedding, implies that the inclusion $F(B_h)\subset B_h$ is compact. Thus, by Schauder's fixed point theorem (c.f. \cite[Theorem 3.1]{Elliott2017}), there exists a fixed point $x_h$ for \eqref{Weak, ODE} that satisfies $F(x_h)=x_h$, which is the desired semi-discrete solution. Moreover, the estimate \eqref{Stability estimate for y_h} also holds for the solution $x_h$.

To address the uniqueness, it is important to recognize that \eqref{Semidiscrete2,weak} constitutes a nonlinear ODE system for $x_j$. Consequently, the uniqueness of $x_h$ is assured by nonlinear ODE theory. It's evident that $x_h$ serves as a semi-discrete solution of \eqref{Semidiscrete2,weak} and thus aligns with the corresponding estimate.
 \qed

\section{Convergence under manifold distance}

\smallskip

As discussed in  \cite{Zhao2021,Jiang2023}, for two closed simple curves $\Gamma_1$ and $\Gamma_2$, the manifold distance is defined as:
\[
\mathrm{M}\l(\Gamma_1,\Gamma_2\r): = \mathrm{Area}((\Omega_1\setminus\Omega_2)\cup (\Omega_2\setminus\Omega_1)) =\mathrm{Area}(\Omega_1)+\mathrm{Area}(\Omega_2 )-2 \mathrm{Area}(\Omega_1\cap \Omega_2 ),
		\]
where $\Omega_1$ and $\Omega_2$ are the  regions enclosed by $\Gamma_1$ and $\Gamma_2$ respectively. As proved in \cite[Proposition 5.1]{Zhao2021}, the manifold distance satisfies symmetry, positivity and  the triangle inequality. Under some assumptions, e.g., If $\Gamma_2$ lies within the tabular neighborhood of the $C^2$ curve $\Gamma_1$ \cite{Deckelnick2005}, the manifold  distance between the two curves can be interpreted as the $L^1$-norm of the distance function. Recently, compare to the $L^p$-norm of parametrization functions, this type of distance (i.e. the $L^p$-norm of distance function) has gained wide attentions in both the scientific computing~\cite{Jiang2023,Jiang2024stable}  and numerical analysis community~\cite{Bai2023}. Moreover, the authors' works \cite{Zhao2021,Jiang2023,Jiang2024stable} have demonstrated that the manifold distance (one of the shape metrics) is  more suitable than the norm of parametrization functions for quantifying numerical errors of the schemes which are used for solving geometric flows, especially for schemes
which allow intrinsic tangential velocity.  Meanwhile, Bai and Li \cite{Bai2023} have recently observed that the $L^2$-norm of distance function (so-called the projected distance in their paper) leads to the recovery of full $H^1$ parabolicity, and established  a convergence result for Dziuk's scheme of the mean curvature flow with finite elements of degree $k\ge 3$.

In this subsection, we first show that the function  $L^\infty$-norm is stronger than the manifold distance under some suitable regularity assumptions. More specifically, for a parametrization function $X\in C^2(\S^1)$ of curve $\Gamma_X$ and  an approximation  curve
$\Gamma_{Y}$ by parameterization function $Y\in C^0(\S^1)$, we have the following lemma.

\begin{lemma}\label{Convergence under mfd, lemma}
Let $X:\S^1\rightarrow \R^2$ be a  parametrization function of simple curves $\Gamma_X$ with $X\in C^{2}(\S^1)$, and assume there exist constants $0<C_1<C_2$ such that \eqref{C12} holds. Then there exist positive constants  $\delta_0$ and $C$ such that for any parametrization function $Y\in C^{0}(\S^1)$ satisfying
\[
	\|X-Y\|_{L^\infty}\le \delta_0,
\]
the following inequality holds:
\[
	\mathrm{M}(\Gamma_X,\Gamma_{Y})\le C\|X-Y\|_{L^\infty},
\]
where $\Gamma_X$ and $\Gamma_{Y}$ are the images of $X$ and $Y$, respectively. The constants $\delta_0$ and $C$ depend on $X$, $C_1$ and $C_2$.
\end{lemma}
\Proof
The closed simple $C^2$ curve $\Gamma_X$ in  $\R^2$ naturally  admits a tabular neighborhood $\Omega_{\delta}$ in the following manner \cite{Bansch2023,Deckelnick2005}: there exists a constant $\delta>0$ such that the mapping
\[
	E_X: \Gamma_{X}\times (-\delta,\delta ) \rightarrow \R^2,\quad
	E_X(a,\eta)=a+\eta\cN,
\]
acts as a diffeomorphism from $\Gamma_{X}\times (-\delta,\delta )$ to the image denoted by $\mathrm{Im}(E_X)=:\Omega_{\delta}$. Here $\cN$ represents the normal vector along $\Gamma_X$. Consequently, the points within the tabular neighborhood $\Omega_\delta$ can be represented as
\[
	E_X^{-1}:\Omega_\delta  \rightarrow \Gamma_{X}\times (-\delta,\delta ),\quad E_X^{-1}(b)=(\pi_{\Gamma_X}(b),d_{\Gamma_X}(b)),
\]
where $\pi_{\Gamma_X}(b)\in \Gamma_X$ is the projection of $b$ onto $\Gamma_X$, and $d_{\Gamma_X}(b)=d(b,\Gamma_X)$ represents the signed distance.

\begin{figure}[htpb]
\hspace{-2mm}
\begin{minipage}[c]{0.8\linewidth}
\centering
\includegraphics[width=14.5cm,height=3.6cm]{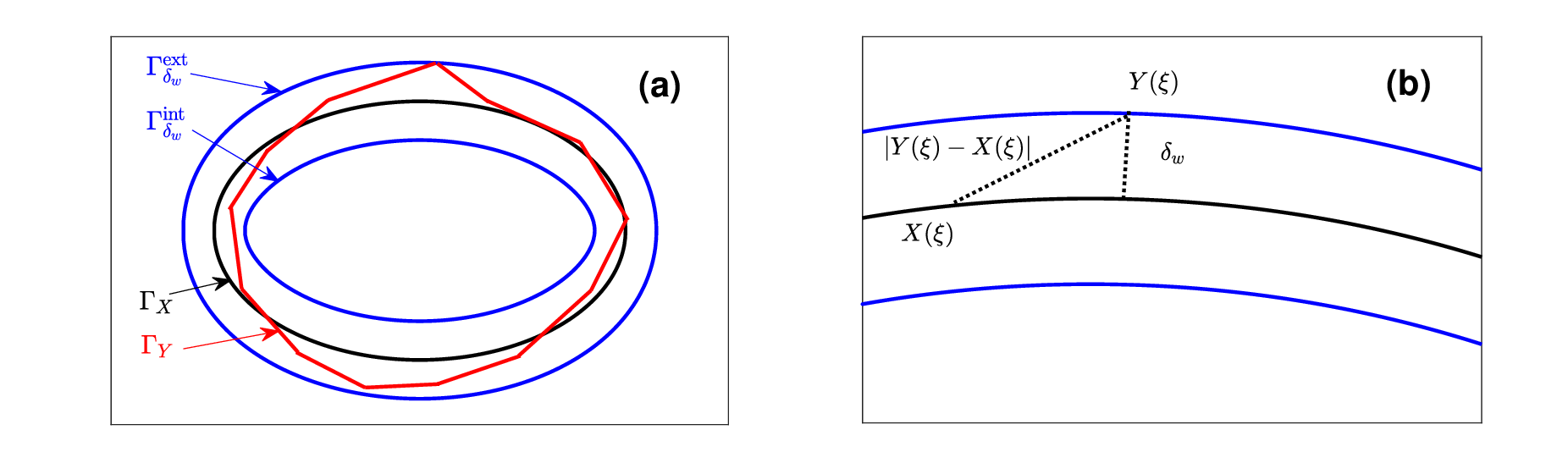}
\setlength{\abovecaptionskip}{-30pt}
\end{minipage}
		\caption{Illustration of (a) the definition of $\Gamma_{X},\Gamma_{Y},\Gamma_{\delta_w}^{\mathrm{int}}$ and $\Gamma_{\delta_w}^\mathrm{{ext}}$, (b) the comparison of the projection distance $\delta_w$ and the function $L^\infty$-norm $\|X-Y\|_{L^\infty}$.}
		\label{Illustration}
	\end{figure}

Set $\delta_0<\delta$. For any parametrization function $Y\in C^0(\S^1)$ which satisfies $\|X-Y\|_{L^\infty}\le \delta_0$, it is evident that $\Gamma_{Y}\subseteq \Omega_{\delta}$. Now define
\[\delta_w:=\sup_{b\in\Gamma_{Y}}|d_{\Gamma_{X}}(b)|\]
which represents the maximum distance between $\Gamma_Y$ and
$\Gamma_X$. Clearly, we can assume $\delta_w>0$, as there is nothing to prove otherwise. Define  two curves $\Gamma_{\delta_w}^{\mathrm{int}}$ and $\Gamma_{\delta_w}^\mathrm{{ext}}$, within the tabular neighborhood $\Omega_{\delta}$, which are  parametrized as
\begin{equation}\label{Parametrization of Gamma boundary}
	\begin{split}
		\S^1\ni \xi \rightarrow (x_{\mathrm{int}}(\xi),y_{\mathrm{int}}(\xi))=(x(\xi),y(\xi))-\delta_w \cN(x(\xi),y(\xi) ),\\
	\S^1\ni \xi \rightarrow (x_{\mathrm{ext}}(\xi),y_{\mathrm{ext}}(\xi))=(x(\xi),y(\xi))+\delta_w \cN(x(\xi),y(\xi) ),
	\end{split}
\end{equation}
respectively. Here $X(\xi)=(x(\xi),y(\xi))$ is a parametrization of the curve $\Gamma_X$, and  $\cN$ is the corresponding unit normal vector. Denote $\Omega_{\delta_w}$ as the region enclosed by $\Gamma_{\delta_w}^{\mathrm{int}}$ and $\Gamma_{\delta_w}^\mathrm{{ext}}$ (cf.  Fig  \ref{Illustration} (a)).
 By utilizing the regularity assumption of $X$ along with \eqref{C12} and \eqref{Parametrization of Gamma boundary}, we can estimate  the area of $\Omega_{\delta_w}$ as follows:
\begin{align*}
	\mathrm{Area}(\Omega_{\delta_w})
	&=\int_{\S^1}\p_\xi x_{\mathrm{ext}}\cdot y_{\mathrm{ext}}\ \d\xi-\int_{\S^1}\p_\xi x_{\mathrm{int}}\cdot y_{\mathrm{int}}\ \d\xi\\
	&\le C\int_{\S^1} |\p_\xi x_{\mathrm{ext}}-\p_\xi x_{\mathrm{int}}|+| y_{\mathrm{ext}}- y_{\mathrm{int}}|\ \d\xi\\
	&\le C\delta_w\int_{\S^1}|\p_\xi\cN |+|\cN| \ \d\xi\le  C\delta_w,
\end{align*}
where $C$ is a constant depending on $X,C_1$ and $C_2$. The triangle inequality for manifold distance yields
\[
	\mathrm{M}(\Gamma_X,\Gamma_{Y})\le 	\mathrm{M}(\Gamma_X,\Gamma_{\delta_w}^{\mathrm{int}})+\mathrm{M}(\Gamma_{\delta_w}^{\mathrm{int}},\Gamma_{Y})\le 2\mathrm{Area}(\Omega_{\delta_w})\le  C\delta_w\le C\|X-Y\|_{L^\infty},\]
where we use the natural control $\delta_{w}\le \|X-Y\|_{L^\infty}$ (cf. Fig \ref{Illustration} (b)) and the proof is completed.\qed

As natural corollaries, we have following convergence results of numerical schemes under the manifold distance.

\begin{corollary}\label{Convergence under manifold distance}
	We have the  following convergence results under the manifold distance.
	\begin{enumerate}
		\item[(1)]  Let $X(\xi,t)$ be the solution of \eqref{Nonlocal,parametrization1}  satisfying Assumption 2.1 and $X\in W^{3,2}(\S^1)$. Let  $x_h(t)$ be the unique finite difference semi-discrete solution of \eqref{Semidiscrete,FDM}. Then we have
	\[
\sup_{t\in [0,T]}\mathrm{M}(\Gamma_X,\Gamma_{x_h})   \leq C h^2.
\]
\item[(2)] Let $X(\xi,t)$ be the solution of \eqref{Nonlocal,parametrization1} satisfying Assumption 4.2 and $X\in W^{3,2}(\S^1)$. Assume that the partition of $\S^{1}$ satisfies Assumption 4.1 and $x_h(t)$ is the unique finite element semi-discrete solution of \eqref{Semidiscrete1,weak}. Then we have
	\[
\sup_{t\in [0,T]}\mathrm{M}(\Gamma_X,\Gamma_{x_h})   \leq C h.
\]

\item[(3)] Let $X(\xi,t)$ be the solution of \eqref{Nonlocal, tangent motion}  satisfying Assumption 4.3 and $X\in W^{3,2}(\S^1)$. Assume that the partition of $\S^{1}$ satisfies Assumption 4.1 and $x_h(t)$ is the unique finite element semi-discrete solution of \eqref{Semidiscrete2,weak}. Then we have
	\[
\sup_{t\in [0,T]}\mathrm{M}(\Gamma_X,\Gamma_{x_h})   \leq C h.
\]
	\end{enumerate}
For all the estimates, $\Gamma_X$ and $\Gamma_{x_h}$ represent the images of $X$ and $x_h$, respectively, and the constant $C$ depends on $C_1,C_2,T$, $f$ and additionally, $K_1(X)$ for (1) and $K_2(X)$ for (2) and (3).	
\end{corollary}

\Proof
~For the first conclusion, combining the Sobolev embedding, triangle inequality, interpolation error, Lemma \ref{Lemma H^1_G} with the main error estimate \eqref{Main estimate 3}, one obtains
\begin{align*}
	\|X-x_h\|_{L^\infty}
	&\le C\|X-x_h\|_{H^1}\le C\|X-I_hX\|_{H^1}+ C\|I_hX-x_h\|_{H^1}\\
	&\le Ch^2\|X\|_{W^{3,2}}+C\|I_hX-x_h\|_{H^1_G}\sqrt{1+h^2/6}\le Ch^2,
\end{align*}
where for the third inequality we have utilized  Lemma \ref{Lemma H^1_G} for the grid function $I_hX-x_h$. Hence, by  applying Lemma \ref{Convergence under mfd, lemma} with $\delta_0=Ch^2$, $Y=x_h$ and {for different time $t\in [0,T]$}, we conclude the first assertion (1). The latter two statements can be similarly confirmed by referring to Theorems \ref{Error estimate 1} and \ref{Error estimate 2}, along with the Sobolev embedding
\[
\sup_{t\in [0,T]}\|X-x_h\|_{L^\infty}\le C\sup_{t\in [0,T]}\|X-x_h\|_{H^1}\le Ch,
\]
and the proof is completed.
\qed

\begin{lemma}\label{Lemma H^1_G}
Let $g:\cG_h\rightarrow \R$ be a grid function. Then we have
\[
	\|g\|_{H^1}^2\le \|g\|_{H^1_G}^2\l(1+h^2/6\r) ,
\]
where we identify a grid function $g$ with the piecewise linear function  over $\S^1$ that connects the grid values of $g$.
\end{lemma}

\Proof
~Denote $M_j=\max\limits_{\xi\in (\xi_j,\xi_{j+1})}\l(|g|^2\r)''$ and  by applying the trapezoidal rule,  we have
\begin{align*}
	\|g\|_{H^1}^2
	&=\sum_{j=1}^N\int_{\xi_j}^{\xi_{j+1}} |g|^2+|\p_\xi g|^2 \ \d\xi\\
	&\le \sum_{j=1}^N\l(\frac{|g(\xi_j)|^2+|g(\xi_{j+1})|^2}{2}h+\frac{h^3}{12}M_j\r)+h\sum_{j=1}^N |\delta g_j|^2\\
	&= h \sum_{j=1}^N \bigl( | g_j |^2 + | \delta g_j |^2 \bigr)+ \frac{h^3}{12}\sum_{j=1}^NM_j.
\end{align*}
Noticing $g$ is a piecewise linear function, we have
 \[
	M_j=\max\limits_{\xi\in (\xi_j,\xi_{j+1})}\l(|g|^2\r)''=
2\max_{\xi\in (\xi_j,\xi_{j+1})} |g'|^2=2|\delta g_{j+1}|^2,
\]
which yields
\[
	\|g\|_{H^1}^2
	\le h \sum_{j=1}^N \bigl( | g_j |^2 + | \delta g_j |^2 \bigr)+ \frac{ h^3}{6}\sum_{j=1}^N|\delta g_{j+1}|^2
	\le \|g\|_{H^1_G}^2\l(1+\frac{1}{6}h^2\r),
\]
and the proof is completed.
 \qed

\section{Numerical results}
In this section, we present numerous numerical experiments for the proposed three different schemes applied to various geometric flows involving the nonlocal term $f(L)$. We first provide full discretizations for the three schemes using backward Euler time discretization. Specifically, we choose an integer $m$, set the time step $\tau=T/m$ and $t_k=k\tau$ for $k=0,\ldots,m$. Given a fixed mesh size $h$ and a time step $\tau=O(h^2)$, we consider the following three cases.
\begin{itemize}
	\item[(i)] For the finite difference method  \eqref{Semidiscrete,FDM}, given $x_h^0=I_hX^0$, for $k\ge 1$, we consider the solution $x_{h}^k\in \cG_h$ of  the following equation ({\textbf{denoted as FDM}})
\begin{equation}\label{Discrete,FDM}
	\frac{x_j^{k}-x_j^{k-1}}{\tau}
     =  \frac{2}{q^{k-1}_j+q_{j+1}^{k-1}}\l(\frac{x^{k}_{j+1}-x^{k}_{j}}{q^{k-1}_{j+1}}-
     \frac{x^{k}_j-x^{k}_{j-1}}{q^{k-1}_j}\r)-f(l_h^{k-1})\frac{n_j^{k-1}+n^{k-1}_{j+1}}{|n^{k-1}_j+n^{k-1}_{j+1}|},
\end{equation}
where $x_j^k$ represents the grid value, $l_h^{k}$ is the perimeter of the polygon with vertices $\{x_j^k\}_j$, and $n_j^{k}=(\tau_j^{k})^\perp$ is the discrete normal vector.

\item[(ii)] For the finite element method \eqref{Semidiscrete1,weak}, given $x_h^0=I_hX^0$, for $k\ge 1$, we consider the solution   $x_h^k=\sum\limits_{j=1}^Nx^k_j\varphi_j\in V_h$ which satisfies ({\textbf{denoted as FEM}})
		\begin{align*}
			&\int_{\S^1}\left|\partial_{\xi} x_{h}^{k-1}\right| \delta_{\tau} x_{h}^{k} \cdot v_{h} \d \xi+\int_{\S^1} \partial_{\xi} x_{h}^{k} \cdot \partial_{\xi} v_{h}/\left|\partial_{\xi} x_{h}^{k-1}\right| \d \xi\\
			&+\int_{\S^1}\mathbf{h}^2|\p_\xi x^{k-1}_h|\p_\xi\delta_{\tau} x^k_h\cdot \p_\xi v_h/6\ \d \xi+\int_{\S^1}f(l_h^{k-1}) (\p_\xi x_h^k)^\perp\cdot v_h \d \xi=0, \quad \forall\ v_{h} \in V_h,
		\end{align*}		
where $\delta_{\tau}$ is the backward finite difference $\delta_{\tau} x_h^k=(x_h^k-x_h^{k-1})/\tau$, and $l_h^{k-1}$ is the length of the image of $x_h^{k-1}$. Or it can be written equivalently as a discretization for the ODE system \eqref{Semidiscrete1,lump mass}:
\begin{equation}\label{Discrete,FEM}
	\frac{q_j^{k-1}+q_{j+1}^{k-1}}{2\tau}(x_j^k-x_j^{k-1})-\frac{x^k_{j+1}-
x^k_{j}}{q^{k-1}_{j+1}}+\frac{x^k_{j}-x^k_{j-1}}{q^{k-1}_{j}}+
\frac{f(l_h^{k-1})}{2}\l(x^k_{j+1}-x^k_{j-1} \r)^\perp=0.
\end{equation}
\item[(iii)] For the finite element method with tangent motions \eqref{Semidiscrete2,weak}, given $x_h^0=I_hX^0$, for fixed $\alpha \in (0,1]$ and $k\ge 1$, $x_h^k=\sum\limits_{j=1}^Nx^k_j\varphi_j\in V_h$ is the solution of the following ({\textbf{denoted as FEM-TM}})
\begin{equation*}
\begin{split}
	&\int_{\S^1}I_h\l[(\alpha\cdot \delta_{\tau} x_h^k +(1-\alpha)\cdot (\delta_{\tau} x_h^k\cdot n_h^{k-1})n_h^{k-1})\cdot  v_h\r]|\p_\xi x_h^{k-1} |^2 \ \d \xi +\int_{\S^1}\p_\xi x_h^k\cdot \p_\xi v_h \ \d \xi \\
	&\qquad \qquad +\int_{\S^1} I_h\l[f(l_h^{k-1})n_h^{k-1}\cdot  v_h\r]|\p_\xi x_h^{k-1}|^2 \ \d \xi =0,\quad \forall \ v_h\in V_h,
\end{split}	
\end{equation*}
where $n_{h}^{k-1}=\l(\frac{\p_\xi x_h^{k-1}}{|\p_\xi x_h^{k-1}|}\r)^{\perp}$ is the unit normal vector. Through a straightforward computation, we find it can be written equivalently as
\begin{equation}\label{Discrete,FEMt}
	\alpha  \frac{x^k_j-x^{k-1}_j}{\tau}+(1-\alpha)\bigg(\frac{x^k_j-x^{k-1}_j}{\tau}\cdot n^{k-1}_j\bigg)n^{k-1}_j =\frac{2(x^k_{j+1}-2x^k_j+x^k_{j-1})}{(q_{j}^{k-1})^2+(q_{j+1}^{k-1})^2}-
 f(l_h^{k-1})n_j^{k-1},
\end{equation}
where $n_j^{k-1}=\l(\frac{x^{k-1}_j-x^{k-1}_{j-1}}{q^{k-1}_j}\r)^\perp.$
\end{itemize}


\subsection{Accuracy test}
To evaluate the convergence order of the proposed three schemes, we primarily consider the following cases of geometric flows with different initial curves:

\noindent\textbf{Case 1}: An ellipse initial curve, parameterized by $(2\cos\theta, \sin\theta)^T$, $\theta\in [0,2\pi]$, with the corresponding flow being the AP-CSF with $f(L)=2\pi/L$;

\noindent\textbf{Case 2}: A four-leaf rose initial curve, parameterized by $(\cos(2\theta)\cos\theta, \cos(2\theta)\sin\theta)^T$, $\theta\in [0,2\pi]$, with
	the corresponding flow being the AP-CSF for nonsimple curves with $f(L)=2\pi\,\mathrm{ind}/L,\mathrm{ind}(\Gamma)=3$;

\noindent\textbf{Case 3}: An ellipse initial curve, with the corresponding flow being a curve flow with area decreasing rate of $\pi$, i.e., $f(L)=(2\pi-\beta)/L,\beta=\pi$.

As the exact solutions of the above cases are unknown, we consider the following numerical errors for the FDM \eqref{Discrete,FDM}:
\begin{alignat*}{2}
               &L^\infty_tH_G^1 \text{ error} \quad &&\left(\mathcal{E}_1\right)_{h,\tau}(T):= \max_{1\le k\le T/\tau} \big\|x^k_{h,\tau}-\widehat{x}^{4k}_{h/2,\tau/4}\big\|_{H^1_G},\\
                & \text{Manifold distance} \,\,\,  &&\left(\mathcal{E}_2\right)_{h,\tau}(T):= \mathrm{M}(\Gamma^{T/\tau}_{h,\tau},\Gamma^{4T/\tau}_{h/2,\tau/4} ) ,
\end{alignat*}
where we view $\widehat{x}^{4k}_{h/2,\tau/4}$ as a grid function over $\cG_h$ with grid values $\{x^{4k}_{2j}\}_{j=1}^N$, and the $L^\infty_G$ norm is defined as $\| g \|_{L^\infty_G}: =\max\limits_{j=1,\ldots,N} |g_j |$. Furthermore, the polygons $\Gamma^{T/\tau}_{h,\tau}$ and $\Gamma^{4T/\tau}_{h/2,\tau/4}$ are the images of $x_{h,\tau}^{T/\tau}$ and $x_{h/2,\tau/4}^{4T/\tau}$, respectively.

	\begin{figure}[htpb]
\hspace{-1mm}
\begin{minipage}[c]{0.8\linewidth}
\centering
\includegraphics[width=14.5cm,height=4.1cm]{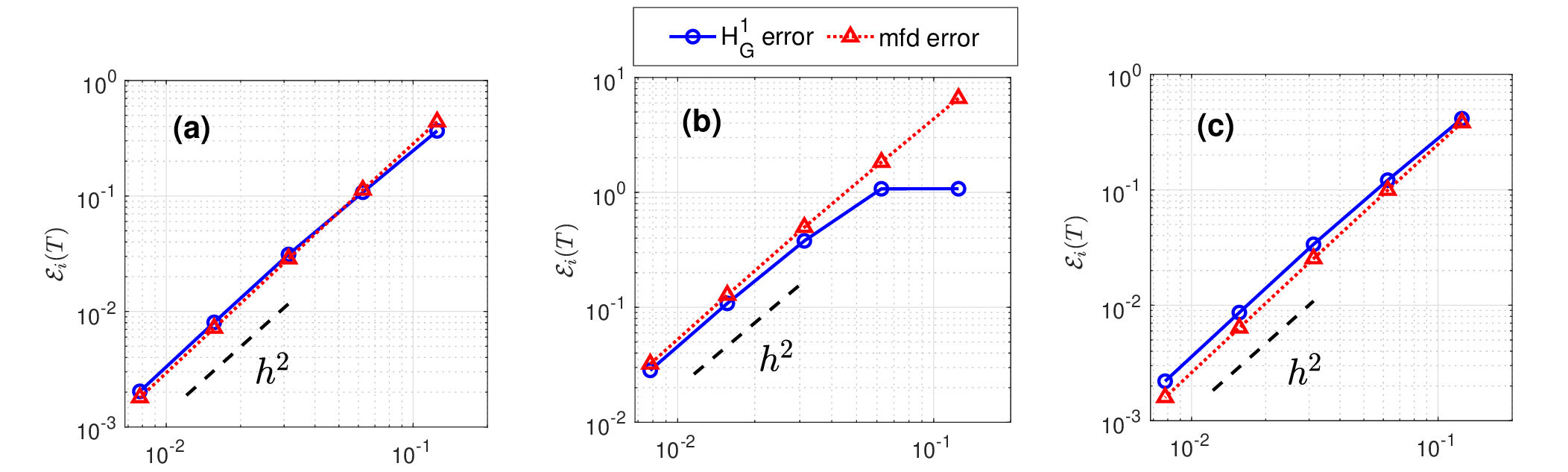}
\setlength{\abovecaptionskip}{-2pt}
\end{minipage}
		\caption{Numerical errors under different norms for the FDM  \eqref{Discrete,FDM} at $T=1/4$: (a) Case 1; (b) Case 2; (c) Case 3.}
		\label{EOC_FDM}
	\end{figure}

Different types of errors for the FDM \eqref{Discrete,FDM} are depicted  Fig \ref{EOC_FDM}, where we choose $h=2\pi/N, \tau=0.5h^2$. The numerical results indicate that, for each instance of nonlocal flows listed above, the solution of \eqref{Discrete,FDM} converges quadratically in $L^\infty_tL^2_G$, $L^\infty_tH^1_G$ and $L^\infty_tL^\infty_G$, which agrees with the theoretical results in Theorem \ref{Error estimate 3}. Moreover, we observe a quadratic convergence under the manifold distance, aligning with the theoretical findings in Corollary \ref{Convergence under manifold distance} (1).

We now turn to the convergence order test of the FEM  \eqref{Discrete,FEM} and the FEM-TM   \eqref{Discrete,FEMt}. We similarly  consider the following numerical errors
\begin{alignat*}{2}
               &L^\infty_tH_x^1 \text{ error} \quad && \left(\mathcal{E}_3\right)_{h,\tau}(T)
                := \max_{1\le k\le T/\tau} \l(\big\|x^k_{h,\tau}-x_{h/2,\tau/4}^{4k}\big\|_{L^2(\mathbb{S}^1)}+\big\|\partial_\xi x^k_{h,\tau}-\partial_\xi x_{h/2,\tau/4}^{4k}\big\|_{L^2(\mathbb{S}^1)}\r),\\
                & \text{Manifold distance}  \,\,\, &&\left(\mathcal{E}_4\right)_{h,\tau}(T):= \mathrm{M}(\Gamma^{T/\tau}_{h,\tau},\Gamma^{4T/\tau}_{h/2,\tau/4} ) ,
\end{alignat*}
where $x_{h,\tau}^k$ represents the solution obtained by the above fully discrete scheme \eqref{Discrete,FEM} or \eqref{Discrete,FEMt} with mesh size $h$ and time step $\tau$.
\begin{figure}[htpb]
\centering
\begin{minipage}[c]{0.8\linewidth}
\hspace{-20mm}
\includegraphics[width=14.5cm,height=4.1cm]{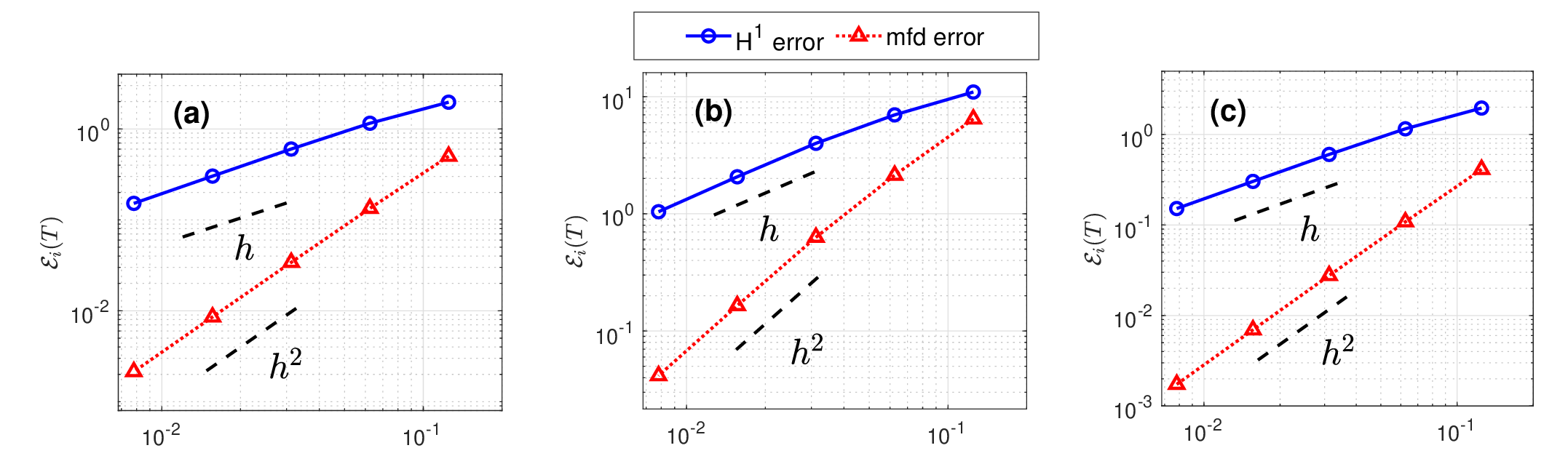}
\setlength{\abovecaptionskip}{-8pt}
\end{minipage}
		\caption{Numerical errors under different norms of the FEM  \eqref{Discrete,FEM} at $T=1/4$: (a) Case 1; (b) Case 2; (c) Case 3.}
		\label{EOC_FEM}
	\end{figure}

\begin{figure}[htpb]
		\hspace{-2mm}
		\includegraphics[width=14.5cm,height=4cm]{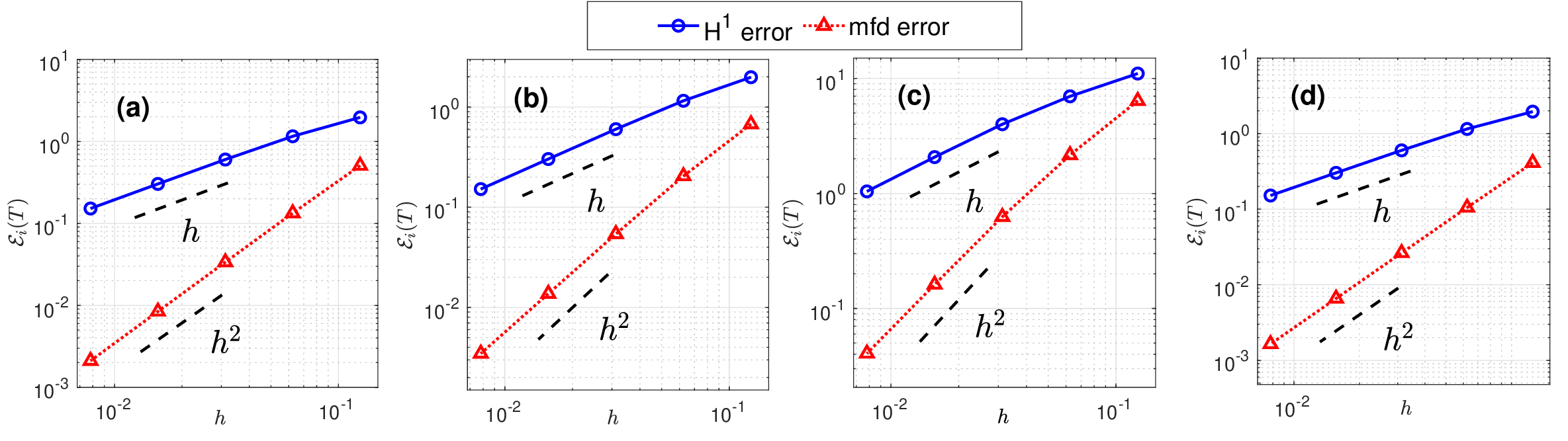}
		\setlength{\abovecaptionskip}{-2pt}
		\caption{Numerical errors under different norms of the FEM-TM  \eqref{Discrete,FEMt} at $T=1/4$: (a) Case 1 with $\alpha=1$; (b) Case 1 with  $\alpha=0.5$; (c) Case 2 with $\alpha=1$; (d) Case 3 with $\alpha=1$.}
		\label{EOC_FEMt}
	\end{figure}

The numerical errors of the FEM \eqref{Discrete,FEM} and the FEM-TM   \eqref{Discrete,FEMt} are presented  in  Fig \ref{EOC_FEM} and Fig \ref{EOC_FEMt}, respectively, from which we observe that, for each nonlocal flow with $\mathrm{ind}$ and  $\beta$, the solution of \eqref{Discrete,FEM} and \eqref{Discrete,FEMt}  converge  linearly in $L^\infty_tH^1_x$, consistent with the theoretical results in Theorems \ref{Error estimate 1} and \ref{Error estimate 2}. Moreover, Fig \ref{EOC_FEMt} (a) and (b) illustrate that the  scheme \eqref{Discrete,FEMt} performs equally well for different choices of $\alpha$.  Additionally, we observe that the solution converges quadratically  under the manifold distance, which is  superior to the theoretical results in Corollary \ref{Convergence under manifold distance} (2) and (3).

\subsection{Dynamics and evolution of geometric quantities}	
In this subsection, we utilize the proposed three  methods: FDM \eqref{Discrete,FDM}, FEM \eqref{Discrete,FEM} and FEM-TM \eqref{Discrete,FEMt} to simulate the nonlocal geometric flows. We are mainly concerned with the evolution of  the following  geometric quantities:  perimeter $L(t)$, relative area loss $\Delta A(t)$ and the mesh ratio function $\Psi(t)$ defined as
\[
\l.L(t)\r|_{t=t_k}=l^k_h,\quad \l.\Delta A(t)\r|_{t=t_k}=\frac{A^k_h-A^0_h}{A^0_h},\quad \l.\Psi(t)\r|_{t=t_k}=\frac{\max_{j=1,\ldots,N}q^k_j}{\min_{j=1,\ldots,N}q^k_j},
\]
where $l^k_h$ and $A^k_h$ are the perimeter and the area  of the polygon determined by $x^k_h$, respectively, and $q^k_j=|x^k_j-x^k_{j-1}|$. Note that for the area of an immersed curve, such as the four-leaf rose, it is treated as a signed area. In morphological evolutions, we primarily focus on the following cases:

\noindent\textbf{Case 1}:  A  flower initial curve parametrized by
\[((2+\cos(6\theta))\cos\theta, (2+\cos(6\theta))\sin\theta)^T,\quad
		 \theta\in [0,2\pi],\]
with the corresponding flow being the AP-CSF with $f(L)=2\pi/L$;

\noindent\textbf{Case 2}:  A four-leaf rose initial curve, with the corresponding flow being the AP-CSF for nonsimple curve with $f(L)=2\pi\,\mathrm{ind}/L,\mathrm{ind}(\Gamma)=3$;

\noindent\textbf{Case 3}: A $4\times 1$ rectangular  initial curve with the corresponding flow being a curve  flow with area decreasing
rate of $\pi$, i.e., $f(L) = (2\pi-\beta)/L$, $\beta=\pi$.

\begin{figure}[htpb]
\hspace{-2mm}
\begin{minipage}[c]{0.8\linewidth}
\centering
\includegraphics[width=14.5cm,height=8cm]{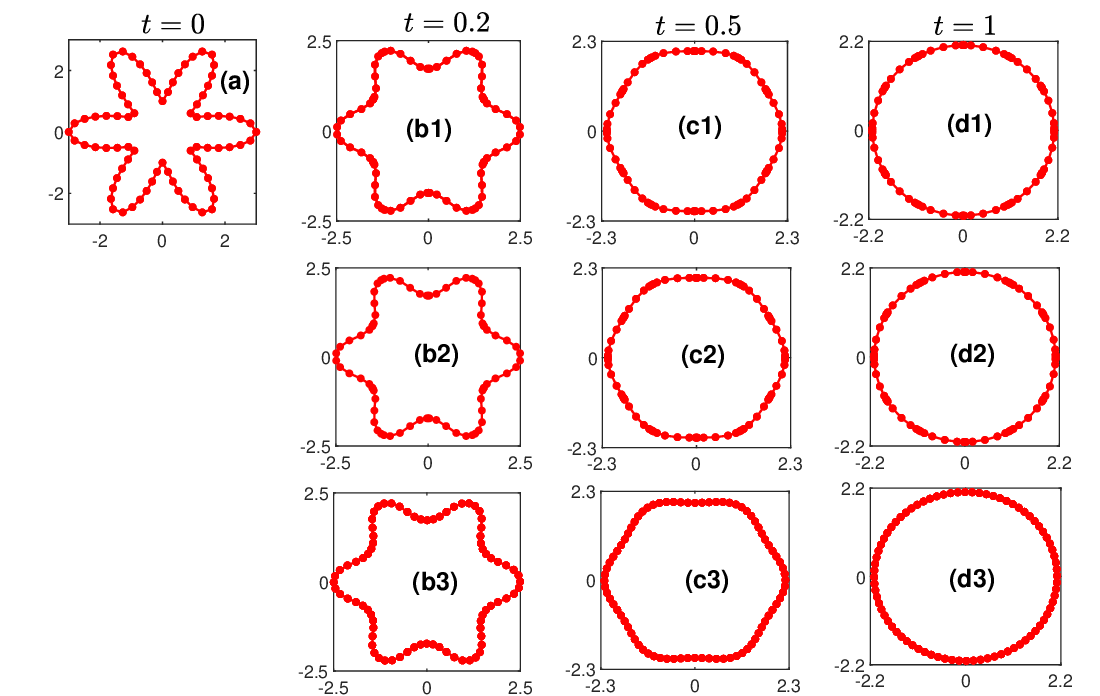}
\setlength{\abovecaptionskip}{-2pt}
\end{minipage}
		\caption{Snapshots of the curve evolution using the FDM (first row), FEM (second row) and FEM-TM (third row) with $\alpha=1$ for  Case 1. The parameters are chosen as $N=80$ and $\tau=1/160$.}
\label{Fig:EVO_case1}		
	\end{figure}

\begin{figure}[htpb]
\hspace{-2mm}
\begin{minipage}[c]{0.8\linewidth}
\centering
\includegraphics[width=14.5cm,height=8cm]{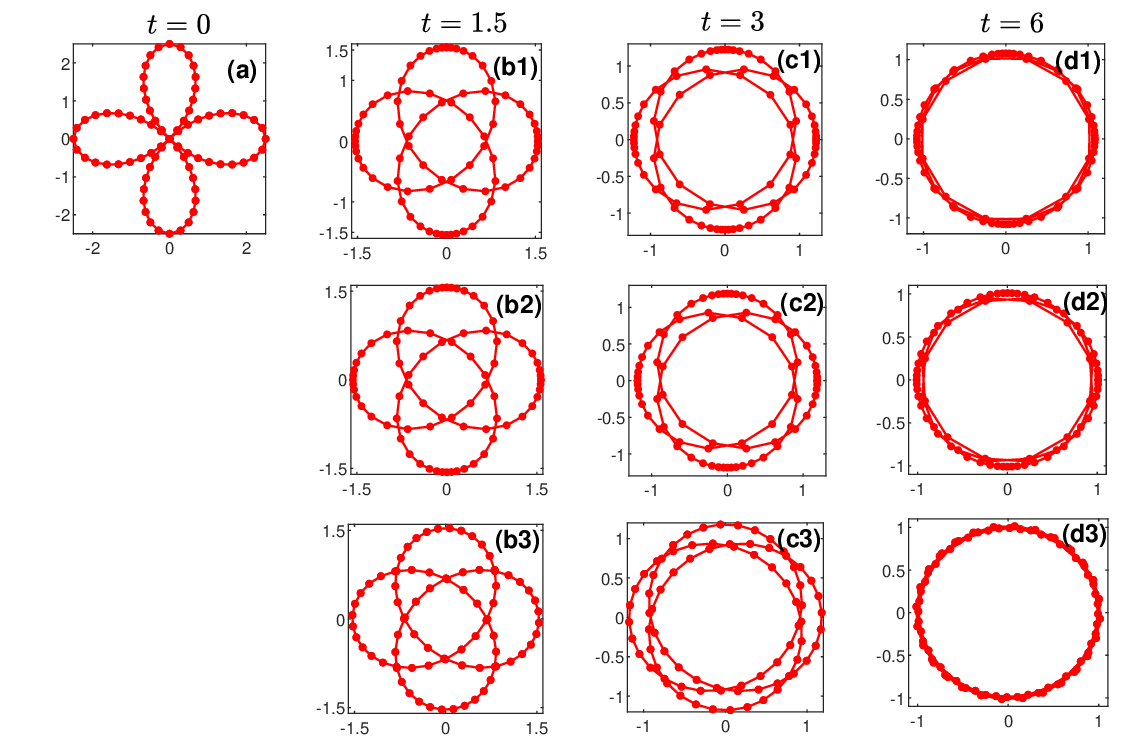}
\setlength{\abovecaptionskip}{-2pt}
\end{minipage}
		\caption{Snapshots of the curve evolution using the FDM (first row), FEM (second row) and FEM-TM (third row) with $\alpha=1$ for  Case 2. The parameters are chosen as $N=80$ and $\tau=1/160$.}
\label{Fig:EVO_case2}
	\end{figure}

\begin{figure}[htpb]
\hspace{-2mm}
\begin{minipage}[c]{0.8\linewidth}
\centering
\includegraphics[width=14.5cm,height=6.5cm]{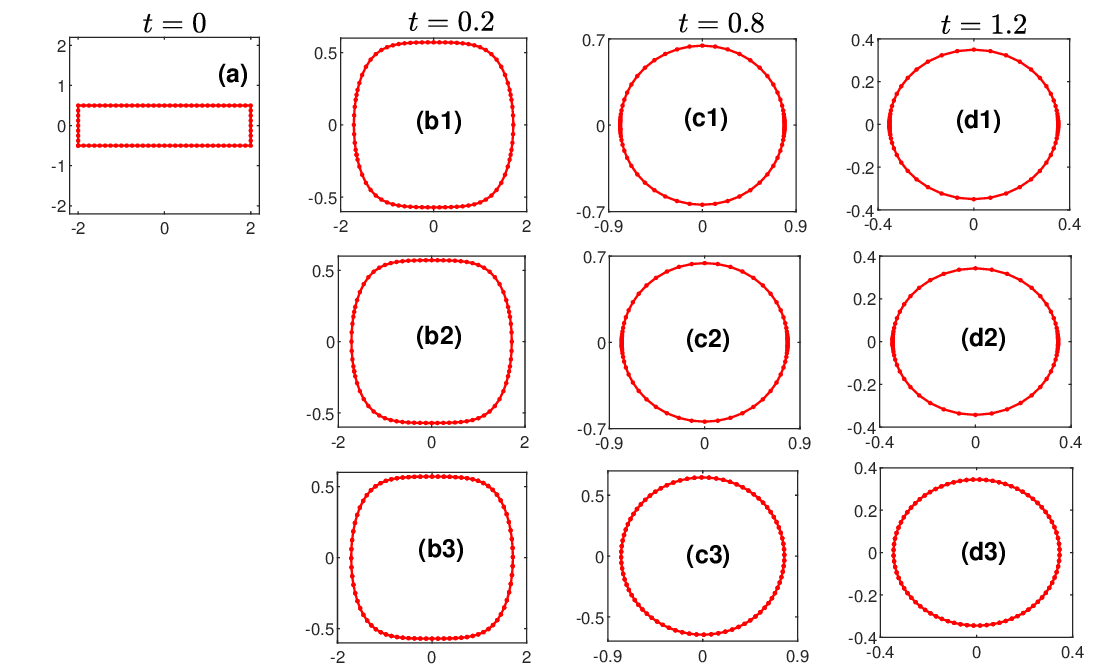}
\setlength{\abovecaptionskip}{-2pt}
\end{minipage}
		\caption{Snapshots of the curve evolution using the FDM (first row), FEM (second row) and FEM-TM (third row) with $\alpha=1$ for  Case 3. The parameters are chosen as $N=80$ and $\tau=1/160$.}
\label{Fig:EVO_case3}		
	\end{figure}

\begin{figure}[htpb]
\hspace{-2mm}
\begin{minipage}[c]{0.8\linewidth}
\centering
\includegraphics[width=14.5cm,height=9.3cm]{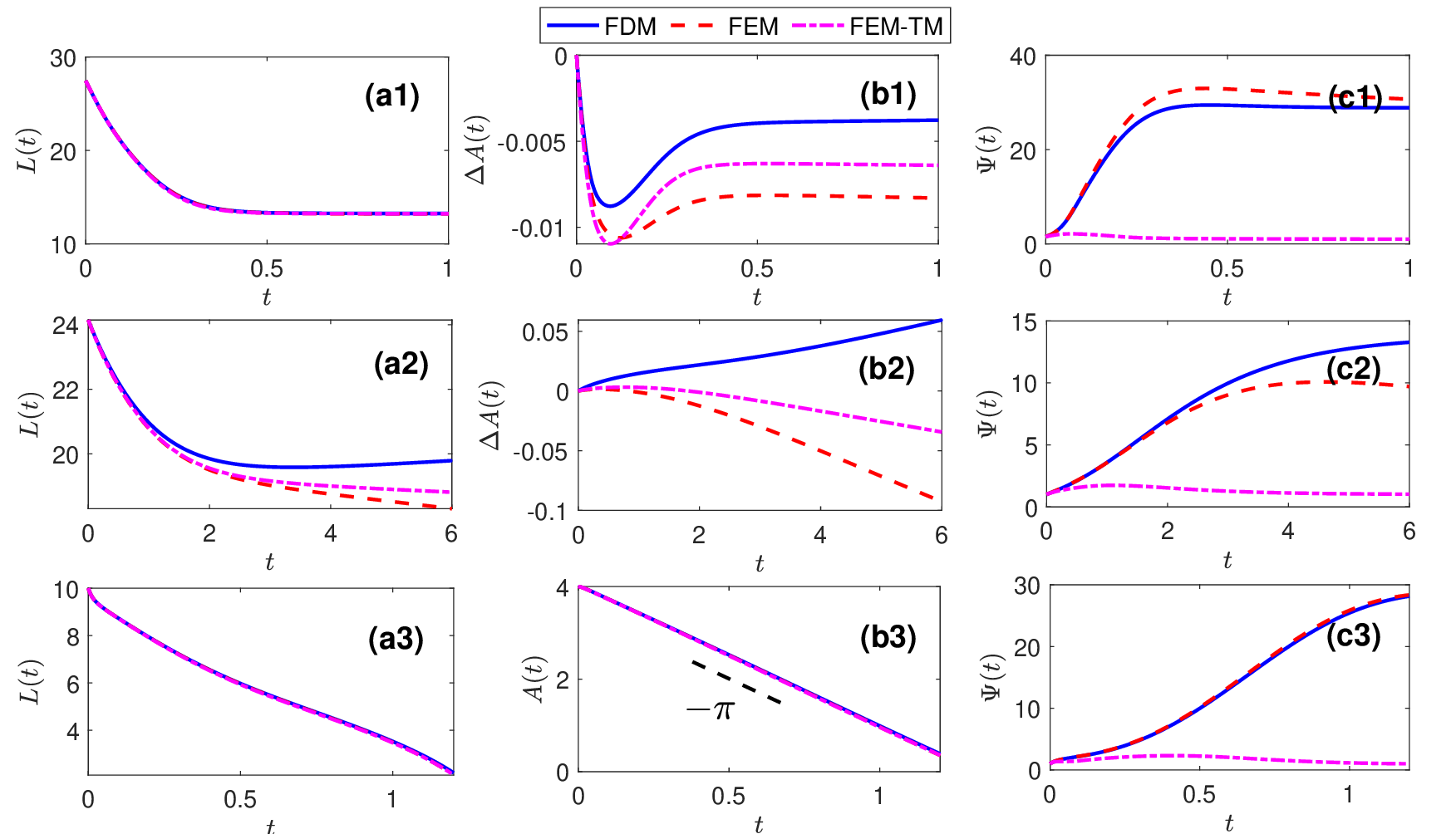}
\end{minipage}
\vspace{4mm}
		\caption{Evolution of the geometric quantities using the FDM, FEM and FEM-TM with $\alpha=1$ for Cases 1-3 is illustrated in the first through third rows, respectively. (a) Perimeter $L(t)$; (b) Relative area loss $\Delta A(t)$; (c) Mesh ratio function $\Psi(t)$. The parameters are chosen as $N=80$ and $\tau=1/160$.}
	\label{Fig:EVO_geo_case123}
	\end{figure}

	Figs. \ref{Fig:EVO_case1}-\ref{Fig:EVO_geo_case123} depict the comparisons of the three schemes through the evolutions of the solution and geometric quantities for the respective three cases. Here we fix the number of nodes $N=80$ and the time step $\tau=1/160$. Based on the observations from Figs. \ref{Fig:EVO_case1}-\ref{Fig:EVO_geo_case123}, we can draw the following conclusions:
\begin{enumerate}
	\item[(i)] All of the schemes can evolve the above three cases successfully into their equilibriums, i.e., circle for Case 1, triple circle for Case 2, and a round point for Case 3, which agrees with the theoretical results in \cite{Wang2014} (cf. Figs. \ref{Fig:EVO_case1}, \ref{Fig:EVO_case2} and \ref{Fig:EVO_case3}).

	\item[(ii)] For Case 1 and Case 2, the area is conserved numerically up to some precision while the area is decreasing numerically with the rate $\pi$ for Case 3 (cf. Fig. \ref{Fig:EVO_geo_case123} (b)).

\item[(iii)] As demonstrated in Fig. \ref{Fig:EVO_geo_case123} (c), the FEM-TM redistributes the points during the evolution  and ultimately achieves the equidistribution property, i.e., $\Psi(t)\rightarrow 1$ as $t\rightarrow +\infty$. This coincides the motivation in Section \ref{sec4}. In contrast, the FDM and the FEM fail to preserve good mesh quality during the evolution.
\end{enumerate}

\begin{figure}[htpb]
\begin{minipage}[c]{0.8\linewidth}
\centering
\includegraphics[width=14cm,height=4.0cm]{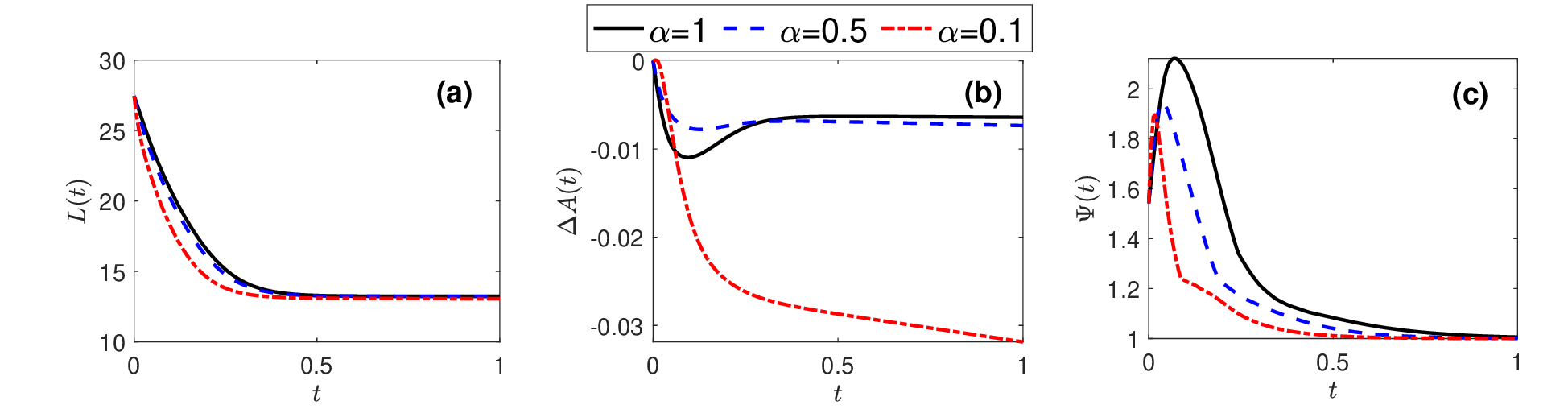}
\end{minipage}
		\caption{Evolution of geometric quantities using the  the FEM-TM with different  $\alpha=0.1,0.5,1$ for  Case 1. (a) Perimeter $L(t)$; (b) Relative area loss $\Delta A(t)$; (c) Mesh ratio function $\Psi(t)$. The parameters are chosen as $N=80$ and $\tau=1/160$.}
	\label{Fig:TMalpha}
	\end{figure}

We close this section with a numerical example to demonstrate that the parameter $\alpha$ in the  FEM-TM \eqref{Discrete,FEMt} signifies the velocity of tangential motions. We conduct simulations for Case 1 using  the FEM-TM with varying values $\alpha=0.1,0.5$ and $1$. As depicted in Fig. \ref{Fig:TMalpha} (c), a smaller $\alpha$
 leads to a more effective redistribution of the mesh points. Fig. \ref{Fig:TMalpha} (b) illustrates that as $\alpha$ approaches $0$, the loss of area becomes greater. This indicates  that for a fixed set of computational parameters $N$ and $\tau$, a smaller value of  $\alpha$ yields a less accurate simulation, aligning with the findings in Theorem \ref{Error estimate 2}, wherein the exponential of $\frac{1}{\alpha}$ is involved in the error estimate.

\section{Conclusions}	

We developed three distinct semi-discrete schemes for simulating some nonlocal geometric flows involving perimeter and the corresponding error estimates were established. Specifically, the FDM exhibits quadratic convergence in $H^1$, whereas the FEM and  the FEM-TM are convergent at the first order in $H^1$. Furthermore, all three methods demonstrate robust quadratic convergence under manifold distance. Extensive numerical experiments have underscored the superior mesh quality of the FEM-TM compared to FDM and FEM.

It is noteworthy that our proof of the error estimate under manifold distance is not optimal for FEM-TM and FEM. Exploring the possibility of providing a  proof of optimal convergence for piecewise linear finite element would be a valuable endeavor.

\section*{Acknowledgments}

W. Jiang was supported by the National Natural Science Foundation of China Grant (No. 12271414) and the Natural Science Foundation of Hubei Province Grant (No. 2022CFB245), and
C. Su and G. Zhang were supported by National Key R\&D Program of China (Grant No. 2023Y-FA1008902) and the National Natural Science Foundation of China Grant (No. 12201342).

\bibliography{reference.bib}

\begin{thebibliography}{10}

\bibitem{Bai2023}
G.~Bai and B.~Li.
\newblock A new approach to the analysis of parametric finite element
  approximations to mean curvature flow.
\newblock {\em Found. Comput. Math., doi.org/10.1007/s10208-023-09622-x}, 2023.

\bibitem{Bansch2023}
E.~B{\"a}nsch, K.~Deckelnick, H.~Garcke, and P.~Pozzi.
\newblock {\em Interfaces: modeling, analysis, numerics. Oberwolfach Seminar,
  Volume 51}.
\newblock Birkh\"auser, Springer, 2023.

\bibitem{Bansch2005}
E.~B\"ansch, P.~Morin, and R.~Nochetto.
\newblock A finite element method for surface diffusion: The parametric case.
\newblock {\em J. Comput. Phys.}, \textbf{203}:321--343, 2005.

\bibitem{Bao2023}
W.~Bao, W.~Jiang, and Y.~Li.
\newblock A symmetrized parametric finite element method for anisotropic
  surface diffusion of closed curves.
\newblock {\em SIAM J. Numer. Anal.}, \textbf{61}(2):617--641., 2023.

\bibitem{Bao2021}
W.~Bao and Q.~Zhao.
\newblock A structure-preserving parametric finite element method for surface
  diffusion.
\newblock {\em SIAM J. Numer. Anal.}, {\textbf{59}}(5):2775--2799., 2021.

\bibitem{Barrett2020}
J.~W. Barrett, H.~Garcke, and R.~N\"urnberg.
\newblock Parametric finite element method approximations of curvature driven
  interface evolutions.
\newblock In Andrea Bonito and Ricardo~H. Nochetto, editors, {\em Handbook of
  Numerical Analysis, Volume 21}, pages 275--423. Elsevier, 2020.

\bibitem{Bronsard1997}
L.~Bronsard and B.~Stoth.
\newblock Volume-preserving mean curvature flow as a limit of a nonlocal
  {G}inzburg-{L}andau equation.
\newblock {\em SIAM J. Math. Anal.}, \textbf{28}(4):769--807., 1997.

\bibitem{Chambolle2015}
A.~Chambolle, M.~Morini, and M.~Ponsiglione.
\newblock Nonlocal curvature flows.
\newblock {\em Arch. Ration. Mech. Anal.}, \textbf{218}:1263--1329., 2015.

\bibitem{Dallaston2012}
M.~C. Dallaston and S.~W. McCue.
\newblock An accurate numerical scheme for the contraction of a bubble in a
  {H}ele-{S}haw cell.
\newblock {\em The ANZIAM Journal}, \textbf{54}:C309--C326, 2012.

\bibitem{Dallaston2013}
M.~C. Dallaston and S.~W. McCue.
\newblock Bubble extinction in {H}ele-{S}haw flow with surface tension and
  kinetic undercooling regularization.
\newblock {\em Nonlinearity}, \textbf{26}(6):1639--1665., 2013.

\bibitem{Dallaston2016}
M.~C. Dallaston and S.~W. McCue.
\newblock A curve shortening flow rule for closed embedded plane curves with a
  prescribed rate of change in enclosed area.
\newblock {\em Proc. R. Soc. A}, \textbf{472}(2185):20150629, 2016.

\bibitem{Deckelnick1995}
K.~Deckelnick and G.~Dziuk.
\newblock On the approximation of the curve shortening flow.
\newblock {\em Pitman Research Notes in Mathematics Series}, pages 100--108.,
  1995.

\bibitem{Deckelnick2005}
K.~Deckelnick, G.~Dziuk, and C.~M. Elliott.
\newblock Computation of geometric partial differential equations and mean
  curvature flow.
\newblock {\em Acta Numer.}, \textbf{14}:139--232., 2005.

\bibitem{Deckelnick2023c}
K.~Deckelnick and R.~N\"urnberg.
\newblock Discrete anisotropic curve shortening flow in higher codimension.
\newblock {\em arXiv:2310.02138}, 2023.

\bibitem{Deckelnick2022}
K.~Deckelnick and R.~N{\"u}rnberg.
\newblock Discrete hyperbolic curvature flow in the plane.
\newblock {\em SIAM J. Numer. Anal.}, \textbf{61}:1835--1857., 2023.

\bibitem{Deckelnick2023b}
K.~Deckelnick and R.~N{\"u}rnberg.
\newblock A novel finite element approximation of anisotropic curve shortening
  flow.
\newblock {\em Interfaces Free Bound.}, \textbf{4}:671--708, 2023.

\bibitem{Deckelnick2023}
K.~Deckelnick and R.~N{\"u}rnberg.
\newblock An unconditionally stable finite element scheme for anisotropic curve
  shortening flow.
\newblock {\em Arch. Math.}, \textbf{59}:263--274., 2023.

\bibitem{Deckelnick2024}
K.~Deckelnick and R.~N\"urnberg.
\newblock Finite element schemes with tangential motion for fourth order
  geometric curve evolutions in arbitrary codimension.
\newblock {\em arXiv:2402.16799}, 2024.

\bibitem{Carmo2016}
do~Carmo M.~P.
\newblock {\em Differential Geometry of Curves and Surfaces}.
\newblock Dover Publications, Inc., Mineola, NY, 2016.

\bibitem{Li2024}
B.~Duan and B.~Li.
\newblock New artificial tangential motions for parametric finite element
  approximation of surface evolution.
\newblock {\em SIAM J. Sci. Comput.}, \textbf{46}:A587--A608, 2024.

\bibitem{Dziuk1999}
G.~Dziuk.
\newblock Discrete anisotropic curve shortening flow.
\newblock {\em SIAM J. Numer. Anal.}, \textbf{36}(6):1808--1830., 1999.

\bibitem{Elliott2017}
C.~M. Elliott and H.~Fritz.
\newblock On approximations of the curve shortening flow and of the mean
  curvature flow based on the {D}e{T}urck trick.
\newblock {\em IMA J. Numer. Anal.}, \textbf{37}(2):543--603., 2017.

\bibitem{Gage1986}
M.~Gage.
\newblock On an area-preserving evolution equation for plane curves.
\newblock {\em Contemp. Math.}, \textbf{51}:51--62., 1986.

\bibitem{Dolcetta2002}
S.~F.~Vita I.~C.~Dolcetta and R.~March.
\newblock Area-preserving curve-shortening flows: from phase separation to
  image processing.
\newblock {\em Interfaces Free Bound.}, \textbf{4}(4):325--343., 2002.

\bibitem{Jiang2008}
L.~Jiang and S.~Pan.
\newblock On a non-local curve evolution problem in the plane.
\newblock {\em Comm. Anal. Geom.}, \textbf{16}(1):1--26., 2008.

\bibitem{Jiang2022}
W.~Jiang, C.~Su, and G.~Zhang.
\newblock A convexity-preserving and perimeter-decreasing parametric finite
  element method for the area-preserving curve shortening flow.
\newblock {\em SIAM J. Numer. Anal.}, \textbf{61}(4):1989--2010., 2023.

\bibitem{Jiang2023}
W.~Jiang, C.~Su, and G.~Zhang.
\newblock A second-order in time, {BGN-based} parametric finite element method
  for geometric flows of curves.
\newblock {\em arXiv:2309.12875.}, 2023.

\bibitem{Jiang2024stable}
W.~Jiang, C.~Su, and G.~Zhang.
\newblock Stable {BDF time discretization of BGN-based} parametric finite
  element methods for geometric flows.
\newblock {\em arXiv:2402.03641}, 2024.

\bibitem{Mayer2000}
U.~F. Mayer.
\newblock A numerical scheme for moving boundary problems that are gradient
  flows for the area functional.
\newblock {\em European J. Appl. Math.}, \textbf{11}(1):61--80., 2000.

\bibitem{Mikula2004b}
K.~Mikula and D.~{\v{S}}ev{\v{c}}ovi{\v{c}}.
\newblock Computational and qualitative aspects of evolution of curves driven
  by curvature and external force.
\newblock {\em Comput. Vis. Sci.}, \textbf{6}(4):211--225., 2004.

\bibitem{Mikula2004}
K.~Mikula and D.~{\v{S}}ev{\v{c}}ovi{\v{c}}.
\newblock A direct method for solving an anisotropic mean curvature flow of
  plane curves with an external force.
\newblock {\em Math. Methods Appl. Sci.}, \textbf{27}(13):1545--1565., 2004.

\bibitem{Pei2022}
L.~Pei and Y.~Li.
\newblock A structure-preserving parametric finite element method for
  area-conserved generalized mean curvature flow.
\newblock {\em J. Sci. Comput.}, \textbf{96}(6):1--21, 2023.

\bibitem{Pozzi2022}
P.~Pozzi and B.~Stinner.
\newblock Convergence of a scheme for elastic flow with tangential mesh
  movement.
\newblock {\em ESAIM Math. Model. Numer. Anal.}, \textbf{57}(2):445--466, 2023.

\bibitem{Ruuth2003}
S.~J. Ruuth and B.~Wetton.
\newblock A simple scheme for volume-preserving motion by mean curvature.
\newblock {\em J. Sci. Comput.}, \textbf{19}:373--384., 2003.

\bibitem{Sapiro2001}
G.~Sapiro.
\newblock {\em Geometric Partial Differential Equations and Image Analysis}.
\newblock Cambridge University Press, 2001.

\bibitem{Sapiro1995}
G.~Sapiro and A.~Tannenbaum.
\newblock Area and length preserving geometric invariant scale-spaces.
\newblock {\em IEEE Trans. Pattern Anal. Mach. Intell.}, \textbf{17}:67--72,
  1995.

\bibitem{Sevcovic2001}
D.~{\v{S}}ev{\v{c}}ovi{\v{c}} and K.~Mikula.
\newblock Evolution of plane curves driven by a nonlinear function of curvature
  and anisotropy.
\newblock {\em SIAM J. Appl. Math.}, \textbf{61}(5):1473--1501., 2001.

\bibitem{Tsai2015}
D.~Tsai and X.~Wang.
\newblock On length-preserving and area-preserving nonlocal flow of convex
  closed plane curves.
\newblock {\em Calc. Var. Partial Differential Equations},
  \textbf{54}:3603--3622., 2015.

\bibitem{Tsai2018}
D.~Tsai and X.~Wang.
\newblock The evolution of nonlocal curvature flow arising in a {H}ele-{S}haw
  problem.
\newblock {\em SIAM J. Math. Anal.}, \textbf{50}(1):1396--1431., 2018.

\bibitem{Ushijima2004}
T.~Ushijima and S.~Yazaki.
\newblock Convergence of a crystalline approximation for an area-preserving
  motion.
\newblock {\em J. Comput. Appl. Math.}, \textbf{166}(2):427--452., 2004.

\bibitem{Wang2014}
X.~Wang and L.~Kong.
\newblock Area-preserving evolution of nonsimple symmetric plane curves.
\newblock {\em J. Evol. Equ.}, \textbf{14}(2):387--401., 2014.

\bibitem{Zhao2021}
Q.~Zhao, W.~Jiang, and W.~Bao.
\newblock An energy-stable parametric finite element method for simulating
  solid-state dewetting.
\newblock {\em IMA J. Numer. Anal.}, \textbf{41}(3):2026--2055., 2021.

\end{thebibliography}

\end{document}